\documentclass[12pt]{article}
\usepackage{amssymb}

\usepackage{graphicx}
\newtheorem{theorem}{{\bf Theorem}}[section]
\newtheorem{lemma}{{\bf Lemma}}[section]
\newtheorem{proposition}{{\bf Proposition}}[section]
\newtheorem{definition}{{\bf Definition}}[section]

\newtheorem{remark}{{\bf Remark}}[section]

\newtheorem{assumption}{{\bf Assumption}}[section]
\begin{document}

\title{Numerical Algorithms and Simulations for Reflected Backward Stochastic
Differential Equations with two Continuous Barriers}
\author{Mingyu XU\thanks{
Email : xumy@amss.ac.cn. This work is supported in part by the
National Basic
Research Program of China (973 Program), No. 2007CB814902} \\
%EndAName
{\small Institute of Applied Mathematics, Academy of Mathematics and Systems
Science,}\\
{\small Chinese Academy of Sciences, Beijing, 100080, China.}}
\date{}
\maketitle

\begin{center}
\textbf{Abstract }
\end{center}

{\small In this paper we study different algorithms for reflected
backward stochastic differential equations (BSDE in short) with two
continuous barriers based on binomial tree framework. We introduce
numerical algorithms by penalization method and reflected method
respectively. In the end
simulation results are also presented. } %\author{}

\vspace*{.8cm} \textbf{Keywords: }{\small Backward Stochastic Differential
Equations with two continuous barriers, Penalization method, Discrete
Brownian motion, Numerical simulation}

\vspace*{.8cm} \textbf{AMS: 60H10, 34K28}

\section{Introduction}

Non-linear backward stochastic differential equations (BSDEs in
short) were firstly introduced by Pardoux and Peng (\cite{PP1990},
1990), who proved the existence and uniqueness of the adapted
solution, under smooth square integrability assumptions on the
coefficient and the terminal condition, plus that the coefficient
$g(t,\omega ,y,z)$ is $(t,\omega )$-uniformly Lipschitz in $(y,z)$.
Then El Karoui, Kapoudjian, Pardoux, Peng and Quenez introduced the
notion of reflected BSDE (RBSDE in short) (\cite{EKPPQ}, 1997) with
one continuous lower barrier. More precisely, a solution for such an
equation associated to a coefficient $g$, a terminal value $\xi $, a
continuous barrier $L_{t}$, is a triplet $(Y_{t},Z_{t},A_{t})_{0\leq
t\leq T} $ of adapted processes valued in $R^{1+d+1}$, which
satisfies
\[
Y_{t}=\xi
+\int_{t}^{T}g(s,Y_{s},Z_{s})ds+A_{T}-A_{t}+\int_{t}^{T}Z_{s}dB_{s},\;\;0%
\leq t\leq T,\mbox{a.s.},
\]
and $Y_{t}\geq L_{t}$ a.s. for any $0\leq t\leq T$. $A_{t}$ is
non-decreasing continuous, and $B_{t}$ is a d-dimensional Brownian
motion. The role of $A_{t}$ is to push upward the process $Y$ in a
minimal way, in
order to keep it above $L$. In this sense it satisfies $%
\int_{0}^{T}(Y_{s}-L_{s})dA_{s}=0$.

Following this paper, Cvitanic and Karatzas (\cite{CK}, 1996)
introduced the notion of reflected BSDE with two continuous
barriers. In this case a solution of such an equation associated to
a coefficient $g$, a terminal value $\xi $, a continuous lower
barrier $L_{t}$ and a continuous upper barrier $U_{t}$, with
$L_{t}\leq U_{t}$ and $L_{T}\leq \xi \leq U_{T}$ a.s., is a
quadruple $(Y_{t},Z_{t},A_{t},K_{t})_{0\leq t\leq T}$ of adapted
processes, valued in $R^{1+d+1}$, which satisfies
\[
Y_{t}=\xi
+\int_{t}^{T}g(s,Y_{s},Z_{s})ds+A_{T}-A_{t}-(K_{T}-K_{t})-%
\int_{t}^{T}Z_{s}dB_{s},0\leq t\leq T,\mbox{a.s.,}
\]
and $L_{t}\leq Y_{t}\leq U_{t}$, a.s. for any $0\leq t\leq T$. Here $A_{t}$
and $K_{t}$ are increasing continuous process, whose roles are to keep the
process $Y$ between $L$ and $U$ in such a way that
\[
\int_{0}^{T}(Y_{s}-L_{s})dA_{s}=0\mbox{ and }%
\int_{0}^{T}(Y_{s}-U_{s})dK_{s}=0.
\]
In view to prove the existence and uniqueness of a solution, the
method is based on a Picard-type iteration procedure, which requires
at each step the solution of a Dynkin game problem. Furthermore, the
authors proved the existence result by penalization method when the
coefficient $g$ does not depend on $z$. In  2004 (\cite{LSM}),
Lepeltier and San Martin relaxed in some sense the condition on the
barriers, proved by a penalization method an existence
result, without any assumption other than the square integrability one on $%
L$ and $U$, but only when there exists a continuous semi-martingale
with terminal value $\xi $, between $L$ and $U$. More recently,
Lepeltier and Xu (\cite{LX04}) studied the case when the barriers
are right continuous and left limit (RCLL in short), and proved the
existence and uniqueness of solution in both Picard iteration and
penalization method. In 2005, Peng and Xu
considered the most general case when barriers are just $\mathbf{L}^{2}$%
-processes by penalization method, and studied a special
penalization BSDE, which penalized with two barriers at the same
time, and proved that the solutions of these equations converge to
the solution of reflected BSDE.

The calculation and simulation of BSDEs is essentially different
from those of SDEs (see \cite{KP}). When $g$ is linear in $y$ and
$z$, we may solve the solution of BSDE by considering its dual
equation, which is a forward SDE. However for nonlinear case of $g$,
we can not find the solution explicitly. Here our numerical
algorithms is based on approximate Brownian motion by random walk.
This method is first considered by Peng and Xu \cite{PX08}. The
convergence of this type of numerical algorithms is proved by
Briand, Delyon and M\'{e}min in 2000 (\cite{bdm}) and 2002
\cite{BDM2} . In 2002, M\'{e}min, Peng and Xu studied the algorithms
for reflected BSDE with one barrier and proved its convergence (cf.
\cite{MPX}). Recently Chassagneux also studied discrete-time
approximation of doubly reflected BSDE in \cite{C09}.

In this paper, we consider different numerical algorithms for
reflected BSDE with two continuous barriers. The basic idea is to
approximate a Brownian motion by random walks based on binary tree
model. Compare with the one barrier case (cf. \cite{MPX}), the
additive barrier brings more difficulties in proving the convergence
of algorithm, which requires us to get finer estimation. When the
Brownian motion is 1-dimensional, our algorithms have advantages in
 computer programming. In fact we developed a software package based on
these algorithms for BSDE with two barriers. Furthermore it also
contains programs for classical BSDEs and reflected BSDEs with one
barrier. One significant advantage of this package is that the users
have a very convenient user-machine interface. Any user who knows
the basics of BSDE can run this package without difficulty.

This paper is organized as follows. In Section 2, we recall some
classical results of reflected BSDE with two continuous barriers,
and discretization for reflected BSDE. In Section 3, we introduce
implicit and implicit-explicit penalization schemes and prove their
convergence. In Section 4, we study implicit and explicit reflected
schemes, and get their convergence. In Section 5, we present some
simulations for reflected BSDE with two barriers. The proof of
convergence of penalization solution is in Appendix.

We should point out that recently there have been many  different
algorithms for computing solutions of BSDEs and the related results
in numerical
analysis, for example \cite{BT}, \cite{bdm}, \cite{C}, \cite{CMM}, \cite{DMP}%
, \cite{GLX}, \cite{MPST}, \cite{ZZ}. In contrast to these results,
our methods can easily be realized by computer in 1-dimensional
case. In the multi-dimensional case, the algorithms are still
suitable, however to realize them by computer is difficult, since it
will require larger amount of calculation than 1-dimensional case.

\section{Preliminaries: Reflected BSDEs with two barriers and Basic
discretization}

Let $(\Omega ,\mathcal{F},P)$ be a complete probability space, $%
(B_{t})_{t\geq 0}$ a $1$-dimensional Brownian motion defined on a fixed
interval $[0,T]$, with a fixed $T>0$. We denote by $\{\mathcal{F}%
_{t}\}_{0\leq t\leq T}$ the natural filtration generated by the Brownian
motion $B$, i.e., $\mathcal{F}_{t}=\sigma \{B_{s};0\leq s\leq t\}$ augmented
with all $P$-null sets of $\mathcal{F}$. Here we mainly consider
1-dimensional case, since the solution of reflected BSDE is 1-dimensional.
In fact, we can also generalize algorithms in this paper to
multi-dimensional Brownian motion, which will require a huge amount of
calculation. We introduce the following spaces for $p\in [1,\infty )$:

\begin{itemize}
\item  $\mathbf{L}^{p}(\mathcal{F}_{t}):=$\{$\mathbb{R}$-valued $\mathcal{F}%
_{t}$--measurable random variables $X$ s. t. $E[|X|^{p}]<\infty $\};

\item  $\mathbf{L}_{\mathcal{F}}^{p}(0,t):=$\{$\mathbb{R}$--valued and $%
\mathcal{F}_{t}$--adapted processes $\varphi $ defined on $[0,t]$, s. t. $%
E\int_{0}^{t}|\varphi _{s}|^{p}ds<\infty $\};

\item  $\mathbf{S}^{p}(0,t):=$\{$\mathbb{R}$--valued and $\mathcal{F}_{t}$%
--adapted continuous processes $\varphi $ defined on $[0,t]$, s. t. $%
E[\sup_{0\leq t\leq T}\left| \varphi _{t}\right| ^{2}]<\infty $\};

\item  $\mathbf{A}^{p}(0,t):=$\{increasing processes in $\mathbf{S}^{p}(0,t)$
with $A(0)=0$\}.
\end{itemize}

\noindent We are especially interested in the case $p=2$.

\subsection{Reflected BSDE: Definition and convergence results}

The random variable $\xi $ is considered as terminal value,
satisfying $\xi
\in \mathbf{L}^{2}(\mathcal{F}_{T})$. Let $g:[0,T]\times \mathbb{R}\times %
\mathbb{R}\rightarrow \mathbb{R}$ be a $t$-uniformly Lipschitz
function in $(y,z)$, i.e., there exists a fixed $\mu >0$ such that
\begin{eqnarray}
|g(t,y_{1},z_{1})-g(t,y_{2},z_{2})| &\leq &\mu (|y_{1}-y_{2}|+|z_{1}-z_{2}|)
\label{Lip} \\
\forall t &\in &[0,T],\forall (y_{1},z_{1}),(y_{2},z_{2})\in \mathbb{R}%
\times \mathbb{R}.  \nonumber
\end{eqnarray}
And $g(\cdot ,0,0)\ $is square integrable.

The solution of our BSDE with two barriers is reflected between a
lower barrier $L$ and an upper barrier $U$, which are supposed to
satisfy

\begin{assumption}
\label{num-bar2} $L$ and $U$ are $\mathcal{F}_{t}$-progressively
measurable continuous processes valued in $\mathbb{R,}$ such that
\begin{equation}
E[\sup_{0\leq t\leq T}((L_{t})^{+})^{2}+\sup_{0\leq t\leq
T}((U_{t})^{-})^{2}]<\infty .  \label{LpUm}
\end{equation}
and there exists a continuous process $X_{t}=X_{0}-\int_{0}^{t}\sigma
_{s}dB_{s}+V_{t}^{+}-V_{t}^{-}$ where $\sigma \in \mathbf{L}_{\mathcal{F}%
}^{2}(0,T)$, $V^{+}$ and $V^{-}$ are $(\mathcal{F}_{t})$-adapted continuous
increasing processes with $E[|V_{T}^{+}|^{2}]+E[|V_{T}^{-}|^{2}]<\infty $
such that
\[
L_{t}\leq X_{t}\leq U_{t},\;\;\;\mbox{ }P\mbox{-a.s. for }0\leq t\leq T.
\]
\end{assumption}

\begin{remark}
Condition (\ref{LpUm}) permits us to treat situations when
$U_{t}\equiv +\infty $ or $L_{t}\equiv $\ $-\infty $, $t\in [0,T]$,
in such cases the corresponding reflected BSDE with two barriers
becomes a reflected BSDE with a single lower barrier $L$ or a single
upper barrier $U$, respectively.
\end{remark}

\begin{definition}
The solution of a reflected BSDE with two continuous barriers is a quadruple $%
(Y,Z,A,K)\in \mathbf{S}^{2}(0,T)\times \mathbf{L}_{\mathcal{F}}^{2}(0,T)$ $%
\times \mathbf{A}^{2}(0,T)\times \mathbf{A}^{2}(0,T)$ defined on
$[0,T]$ satisfying the following equations
\begin{eqnarray}
-dY_{t} &=&g(t,Y_{t},Z_{t})dt+dA_{t}-dK_{t}-Z_{t}dB_{t},\;Y_{T}=\xi
\label{RBSDE2b1} \\
L_{t} &\leq &Y_{t}\leq U_{t},\;\;\;dA_{t}\geq 0,dK_{t}\geq 0,\;\;dA_{t}\cdot
dK_{t}=0.  \nonumber
\end{eqnarray}
and the reflecting conditions
\begin{equation}
\int_{0}^{T}(Y_{t}-L_{t})dA_{t}=\int_{0}^{T}(Y_{t}-U_{t})dK_{t}=0.
\label{refle2}
\end{equation}
\end{definition}

To prove the existence of the solution, penalization method is
important. Thanks to the convergence results of penalization
solution in \cite{LSM}, \cite {HLM} for continuous barriers' case
and methods in \cite{PX2003}, we have the following results,
especially it gives the convergence speed of penalization solutions.

\begin{theorem}
\label{p2b}(a) There exists a unique solution $(Y,Z,A,K)$ of
reflected BSDE, i.e. it satisfies (\ref{RBSDE2b1}), (\ref{refle2}).
Moreover it is the limit of penalization solutions
$(\widehat{Y}_{t}^{m,p},\widehat{Z}_{t}^{m,p},\widehat{A}_{t}^{m,p},\widehat{K}_{t}^{m,p})$
as $m\rightarrow \infty $ then $p\rightarrow \infty $, or equivalent as $%
q\rightarrow \infty $ then $m\rightarrow \infty $. Here the
penalization solution
$(\widehat{Y}_{t}^{m,p},\widehat{Z}_{t}^{m,p},\widehat{A}_{t}^{m,p},\widehat{K}_{t}^{m,p})$
with respect to
two barriers $L$ and $U$ is defined, for $m\in \mathbb{N}$, $p\in \mathbb{N}$%
, as the solution of a classical BSDE
\begin{eqnarray}
-d\widehat{Y}_{t}^{m,p}
&=&g(t,\widehat{Y}_{t}^{m,p},\widehat{Z}_{t}^{m,p})dt+m(\widehat{Y}_{t}^{m,p}-L_{t})^{-}dt-q(\widehat{Y}_{t}^{m,p}-U_{t})^{+}dt-\widehat{Z}_{t}^{m,p}dB_{t},
\label{PBSDE2b-g} \\
\widehat{Y}_{T}^{m,p} &=&\xi .  \nonumber
\end{eqnarray}
And we set $\widehat{A}_{t}^{m,p}=m\int_{0}^{t}(\widehat{Y}_{s}^{m,p}-L_{s})^{-}ds$, $%
\widehat{K}_{t}^{m,p}=p\int_{0}^{t}(\widehat{Y}_{s}^{m,p}-U_{s})^{+}ds$. \\[0.2cm]
(b) Consider a special penalized BSDE for the reflected BSDE with
two barriers: for any $p\in \mathbb{N}$,
\begin{eqnarray}
-dY_{t}^{p}
&=&g(t,Y_{t}^{p},Z_{t}^{p})dt+p(Y_{t}^{p}-L_{t})^{-}dt-p(Y_{t}^{p}-U_{t})^{+}dt-Z_{t}^{p}dB_{t},
\label{PBSDE2b} \\
Y_{T}^{p} &=&\xi ,  \nonumber
\end{eqnarray}
with $A_{t}^{p}=\int_{0}^{t}p(Y_{s}^{p}-L_{s})^{-}ds$ and $%
K_{t}^{p}=\int_{0}^{t}p(Y_{s}^{p}-U_{s})^{+}ds$. Then we have, as $%
p\rightarrow \infty $, $Y_{t}^{p}\rightarrow Y_{t}$ in $\mathbf{S}^{2}(0,T)$%
, $Z_{t}^{p}\rightarrow Z_{t}$ in $\mathbf{L}_{\mathcal{F}}^{2}(0,T)$ and $%
A_{t}^{p}\rightarrow A_{t}$ weakly in $\mathbf{S}^{2}(0,T)$ as well as $%
K_{t}^{p}\rightarrow K_{t}$. Moreover there exists a constant $C$ depending
on $\xi $, $g(t,0,0)$, $\mu $, $L$ and $U$, such that
\begin{equation}
E[\sup_{0\leq t\leq
T}|Y_{t}^{p}-Y_{t}|^{2}+\int_{0}^{T}(|Z_{t}^{p}-Z_{t}|^{2})dt+\sup_{0\leq
t\leq T}[(A_{t}-K_{t})-(A_{t}^{p}-K_{t}^{p})]^{2}]\leq \frac{C}{\sqrt{p}}.
\label{conl2-p}
\end{equation}
\end{theorem}

The proof is based on the results in \cite{LSM} and \cite{PX2003},
we put it in Appendix.
\begin{remark}
In the following, we focus on the penalized BSDE as (\ref{conl2-p}),
which consider the penalization with respect to the two barriers at
the same time. And $p$ in superscribe always stands for the
penalization parameter.
\end{remark}

Now we consider a special case: Assume that

\begin{assumption}
\label{Ito} $L$ and $U$ are It\^{o} processes of the following form
\begin{eqnarray}
L_{t} &=&L_{0}+\int_{0}^{t}l_{s}ds+\int_{0}^{t}\sigma _{s}^{l}dB_{s},
\label{barriers} \\
U_{t} &=&U_{0}+\int_{0}^{t}u_{s}ds+\int_{0}^{t}\sigma _{s}^{u}dB_{s}.
\nonumber
\end{eqnarray}
Suppose that $l_{s}$ and $u_{s}$ are right continuous with left
limits (RCLL in short) processes, $\sigma _{s}^{l}$ and $\sigma
_{s}^{u}$ are predictable with $E[\int_{0}^{T}[\left| l_{s}\right|
^{2}+\left| \sigma _{s}^{l}\right| ^{2}+\left| u_{s}\right|
^{2}+\left| \sigma _{s}^{u}\right| ^{2}]ds<\infty $.
\end{assumption}

It is easy to check that if $L_{t}\leq U_{t}$, then Assumption \ref{num-bar2}
is satisfied. We may just set $X=L$ or $U$, with $\sigma _{s}=\sigma
_{s}^{l} $ or $\sigma _{s}^{u}$ and $V^{\pm }=\int_{0}l_{s}^{\pm }ds$ or $%
\int_{0}u_{s}^{\pm }ds$. Here $l_{s}^{\pm }$ (resp. $u_{s}^{\pm }$)
is the positive or the negative part of $l$ (resp. $u$). As
Proposition 4.2 in \cite {EKPPQ}, we have following proposition for
two increasing processes, which can give us the integrability of the
increasing processes by barriers.

\begin{proposition}
Let $(Y,Z,A,K)$ be a solution of reflected BSDE (\ref{RBSDE2b1}). Then $%
Z_{t}=\sigma _{t}^{l}$, a.s.-$dP\times dt$ on the set $\{Y_{t}=L_{t}\}$, $%
Z_{t}=\sigma _{t}^{u}$, a.s.-$dP\times dt$ on the set
$\{Y_{t}=U_{t}\}$. And
\begin{eqnarray*}
0 &\leq &dA_{t}\leq 1_{\{Y_{t}=L_{t}\}}[g(t,L_{t},\sigma
_{t}^{l})+l_{t}]^{-}dt, \\
0 &\leq &dK_{t}\leq 1_{\{Y_{t}=U_{t}\}}[g(t,U_{t},\sigma
_{t}^{u})+u_{t}]^{+}dt.
\end{eqnarray*}
So there exist positive processes $\alpha $ and $\beta $, with
$0\leq \alpha _{t},\beta _{t}\leq 1$, such that $dA_{t}=\alpha
_{t}1_{\{Y_{t}=L_{t}\}}[g(t,L_{t},\sigma _{t}^{l})+l_{t}]^{-}dt$, $%
dK_{t}=\beta _{t}1_{\{Y_{t}=U_{t}\}}[g(t,U_{t},\sigma
_{t}^{u})+u_{t}]^{+}dt. $
\end{proposition}

\textbf{Proof.} We can prove these results easily by using similar
techniques as in Proposition 4.2 in \cite{EKPPQ}, in view that on
the set $\{L_{t}=U_{t}\}$, we have $\sigma _{t}^{l}=\sigma _{t}^{u}$
and $l_{t}=u_{t}$. So we omit the details of the proof here.
$\square $

In the following, we will work under Assumption \ref{Ito}.

\subsection{Approximation of Brownian motion and barriers}

We use random walk to approximate the Brownian motion. Consider for each $%
j=1,2,\cdots ,$
\[
B_{t}^{n}:=\sqrt{\delta }\sum_{j=1}^{[t/\delta ]}\varepsilon _{j}^{n},\;\;%
\mbox{ for all }\;0\leq t\leq T,\;\delta =\frac{T}{n},
\]
where $\{\varepsilon _{j}^{n}\}_{j=1}^{n}$ is a $\{1,-1\}$-valued
i.i.d. sequence with $P(\varepsilon _{j}^{n}=1)=P(\varepsilon
_{j}^{n}=-1)=0.5$,
i.e., it is a Bernoulli sequence. We set the discrete filtration $\mathcal{G}%
_{j}^{n}:=\sigma \{\varepsilon _{1}^{n},\cdots ,\varepsilon _{j}^{n}\}$ and $%
t_{j}=j\delta $, for $0\leq j\leq n$. We denote by $\mathbf{D}_{t}$
the space of RCLL functions from $[0,t]$ in $\mathbb{R}$, endowed
with the topology of uniform convergence, and we assume that:

\begin{assumption}
\label{ass-ter} $\Gamma :\mathbf{D}_{T}\mathbf{\rightarrow R}$ is $K$%
-Lipschitz. We consider $\xi =\Gamma (B)$, which is $\mathcal{F}_{T}$%
-measurable and $\xi ^{n}=\Gamma (B^{n})$, which is $\mathcal{G}_{n}^{n}$%
-measurable, such that
\[
E[|\xi |^{2}]+\sup_{n}E[|\xi ^{n}|^{2}]<\infty
\]
\end{assumption}

Now we consider the approximation of the barriers $L$ and $U$. Notice that $%
L $ and $U$ are progressively measurable with respect to the filtration $(%
\mathcal{F}_{t})$, which is generated by Brownian motion. So they can be
presented as a functional of Brownian motion, i.e. for each $t\in [0,T]$, $%
L_{t}=\Psi _{1}(t,(B_{s})_{0\leq s\leq t})$ and $U_{t}=\Psi
_{2}(t,(B_{s})_{0\leq s\leq t})$, where $\Psi _{1}(t,\cdot )$ and $\Psi
_{2}(t,\cdot ):\mathbf{D}_{t}\mathbf{\rightarrow R}$. And we assume that $%
\Psi _{1}(t,\cdot )$ and $\Psi _{2}(t,\cdot )$ are Lipschitz. Then
we get the discretizaton of the barriers $L_{j}^{n}=\Psi
_{1}(t_{j},(B_{s}^{n})_{0\leq s\leq t})$ and $U_{t}^{n}=\Psi
_{2}(t_{j},(B_{s}^{n})_{0\leq s\leq t})$. If $L_t\leq U_t$, then
$L_j\leq U_j$. On the other hand, we mainly consider barriers which
are It\^o processes and satisfy Assumption \ref{Ito}. So we have a
natural approximation: for $1\leq j\leq n$,
\begin{eqnarray*}
L_{j}^{n} &=&L_{0}+\delta \sum_{i=0}^{j-1}l_{i}+\sum_{i=0}^{j-1}\sigma
_{i}^{l}\varepsilon _{i+1}^{n}\sqrt{\delta }, \\
U_{j}^{n} &=&U_{0}+\delta
\sum_{i=0}^{j-1}u_{i}+\sum_{i=0}^{j-1}\sigma _{i}^{u}\varepsilon
_{i+1}^{n}\sqrt{\delta }
\end{eqnarray*}
where $l_{i}=l_{t_{i}}$, $\sigma _{i}^{l}=\sigma
_{t_{i}}^{l}$, $u_{i}=u_{t_{i}}$%
, $\sigma _{i}^{u}=\sigma _{t_{i}}^{u}$. Then $%
L_{j}^{n}$ and $U_{j}^{n}$ are discrete versions of $L$ and $U$, with $%
\sup_{n}E[\sup_{j}((L_{j}^{n})^{+})^{2}+\sup_{j}((U_{j}^{n})^{-})^{2}]<%
\infty $ and $L_{j}^{n}\leq U_{j}^{n}$ still hold. In the following,
we may use both approximations.

In this paper, we study two different types of numerical schemes.
The first one is based on the penalization approach, whereas the
second is to obtain the solution $Y$ by reflecting it between $L$
and $U$ and get two reflecting processes $A$ and $K$ directly.
Throughout this paper, $n$ always stands for the discretization of
the time interval. And process $(\phi^n_j)_{0\leq j\leq n}$ is a
discrete process with $n+1$ values, for $\phi = L$, $U$, $y^p$, $y$,
etc.

\section{Algorithms based on Penalization BSDE and their Convergence}

\subsection{Discretization of Penalization BSDE and Penalization schemes}

First we consider the discretization of penalized BSDE with respect to two
discrete barriers $L^{n}$ and $U^{n}$. After the discretization of time
interval, we get the following discrete backward equation on the same interval $%
[t_{j},t_{j+1}]$, for $0\leq j\leq n-1$,
\begin{eqnarray}
y_{j}^{p,n} &=&y_{j+1}^{p,n}+g(t_{j},y_{j}^{p,n},z_{j}^{p,n})\delta
+a_{j}^{p,n}-k_{j}^{p,n}-z_{j}^{p,n}\sqrt{\delta }\varepsilon _{j+1}^{n},
\label{dis-p2b} \\
a_{j}^{p,n} &=&p\delta (y_{j}^{p,n}-L_{j}^{n})^{-},\;\;k_{j}^{p,n}=p\delta
(y_{j}^{p,n}-U_{j}^{n})^{+}.  \nonumber
\end{eqnarray}
The terminal condition is $y_{n}^{p,n}=\xi ^{n}$. Since for a large fixed $%
p>0$, (\ref{PBSDE2b}) is in fact a classical BSDE. By numerical algorithms
for BSDEs (cf. \cite{PX2007}), explicit scheme gives $z_{j}^{p,n}=\frac{1}{2%
\sqrt{\delta }}(y_{j+1}^{n}|_{\varepsilon
_{j}=1}-y_{j+1}^{n}|_{\varepsilon _{j}=-1})$, and $y_{j}^{p,n}$ is
solved from the inversion of the following mapping
\begin{eqnarray*}
&&y_{j}^{p,n}=(\Theta ^{p})^{-1}(E[y_{j+1}^{p,n}|\mathcal{G}_{j}^{n}]), \\
&&\mbox{where }\Theta ^{p}(y)=y-g(t_{j},y,z_{j}^{p,n})\delta -p\delta
(y-L_{j}^{n})^{-}+p\delta (y-U_{j}^{n})^{+},
\end{eqnarray*}
by substituting $E[y_{j+1}^{p,n}|\mathcal{G}_{j}^{n}]=\frac{1}{2}%
(y_{j+1}^{p,n}|_{\varepsilon _{j+1}^{n}=1}+y_{j+1}^{p,n}|_{\varepsilon
_{j+1}^{n}=-1})$ into it. And increasing processes $a_{j}^{p,n}$ and $%
k_{j}^{p,n}$ will be obtained from (\ref{dis-p2b}).

In many cases, the inversion of the operator $\Theta ^{p}$ is not
easy to solve. So we apply the implicit--explicit penalization
scheme to (\ref{dis-p2b}), replacing $y_{j}^{p,n}$ in $g$ by
$E[y_{j+1}^{p,n}|\mathcal{G}_{j}^{n}]$, and get
\begin{eqnarray*}
\bar{y}_{j}^{p,n} &=&\bar{y}_{j+1}^{p,n}+g(t_{j},E[\bar{y}_{j+1}^{p,n}|%
\mathcal{G}_{j}^{n}],\bar{z}_{j}^{p,n})\delta +\overline{a}_{j}^{p,n}-%
\overline{k}_{j}^{p,n}-\bar{z}_{j}^{p,n}\sqrt{\delta }\varepsilon^n_{j+1} \\
\ \overline{a}_{j}^{p,n} &=&p\delta (\bar{y}_{j}^{p,n}-L_{j}^{n})^{-},\;\;%
\overline{k}_{j}^{p,n}=p\delta (\bar{y}_{j}^{p,n}-U_{j}^{n})^{+}.
\end{eqnarray*}
In the same way, we get $\bar{z}_{j}^{p,n}=\frac{1}{2\sqrt{\delta }}(\bar{y}%
_{j+1}^{n}|_{\varepsilon _{j}^{n}=1}-\bar{y}_{j+1}^{n}|_{\varepsilon
_{j}^{n}=-1})$ and
\begin{equation}
\bar{y}_{j}^{p,n}=E[\bar{y}_{j+1}^{p,n}|\mathcal{G}_{j}^{n}]+g(t_{j},E[\bar{y%
}_{j+1}^{p,n}|\mathcal{G}_{j}^{n}],\bar{z}_{j}^{p,n})\delta +\overline{a}%
_{j}^{p,n}-\overline{k}_{j}^{p,n}.  \label{p-expli}
\end{equation}
Solving this equation, we obtain
\begin{eqnarray*}
\overline{y}_{j}^{p,n} &=&E[\overline{y}_{j+1}^{p,n}|\mathcal{G}%
_{j}^{n}]+g(t_{j},E[\overline{y}_{j+1}^{p,n}|\mathcal{G}_{j}^{n}],\overline{z%
}_{j}^{p,n})\delta \\
&&+\frac{p\delta }{1+p\delta }(E[\overline{y}_{j+1}^{p,n}|\mathcal{G}%
_{j}^{n}]+g(t_{j},E[\overline{y}_{j+1}^{p,n}|\mathcal{G}_{j}^{n}],\overline{z%
}_{j}^{p,n})\delta -L_{j}^{n})^{-} \\
&&-\frac{p\delta }{1+p\delta }(E[\overline{y}_{j+1}^{p,n}|\mathcal{G}%
_{j}^{n}]+g(t_{j},E[\overline{y}_{j+1}^{p,n}|\mathcal{G}_{j}^{n}],\overline{z%
}_{j}^{p,n})\delta -U_{j}^{n})^{+}.
\end{eqnarray*}
with $E[\bar{y}_{j+1}^{p,n}|\mathcal{G}_{j}^{n}]=\frac{1}{2}(\bar{y}%
_{j+1}^{p,n}|_{\varepsilon _{j+1}^{n}=1}+\bar{y}_{j+1}^{p,n}|_{\varepsilon
_{j+1}^{n}=-1})$. For increasing processes, we can get them from
\begin{eqnarray*}
\overline{a}_{j}^{p,n} &=&\frac{p\delta }{1+p\delta }(E[\bar{y}_{j+1}^{p,n}|%
\mathcal{G}_{j}^{n}]+g(t_{j},E[\bar{y}_{j+1}^{p,n}|\mathcal{G}_{j}^{n}],\bar{%
z}_{j}^{p,n})\delta -L_{j}^{n})^{-}, \\
\overline{k}_{j}^{p,n} &=&\frac{p\delta }{1+p\delta }(E[\bar{y}_{j+1}^{p,n}|%
\mathcal{G}_{j}^{n}]+g(t_{j},E[\bar{y}_{j+1}^{p,n}|\mathcal{G}_{j}^{n}],\bar{%
z}_{j}^{p,n})\delta -U_{j}^{n})^{+}.
\end{eqnarray*}

\subsection{Convergence of penalization schemes and estimations}

First we give the following lemma, which is proved in \cite{MPX}.
This Gronwall type lemma is classical but here it is given with more
detailed formulation.

\begin{lemma}
\label{dis-gro}Let $a$, $b$ and $\alpha $ be positive constants, $\delta b<1$
and a sequence $(v_{j})_{j=1,\ldots n}$ of positive numbers such that, for
every $j$%
\begin{equation}
v_{j}+\alpha \leq a+b\delta \sum_{i=1}^{j}v_{i}.  \label{dis-gronwall}
\end{equation}
Then
\[
\sup_{j\leq n}v_{j}+\alpha \leq a\mathcal{E}_{\delta }(b),
\]
where $\mathcal{E}_{\delta }(b)=1+\sum_{p=1}^{\infty }\frac{b^{p}}{p}%
(1+\delta )\cdots (1+(p-1)\delta )$, which is a convergent series.
\end{lemma}

Notice the $\mathcal{E}_{\delta }(b)$ is increasing in $\delta $ and $%
\delta <\frac{1}{b}$, so we can replace the right hand side of (\ref{dis-gronwall}%
) by a constant depending on $b$.

We define the discrete solutions,
$(Y_{t}^{p,n},Z_{t}^{p,n},A_{t}^{p,n},\ K_{t}^{p,n})$ by the
implicit penalization scheme
\[
Y_{t}^{p,n}=y_{[t/\delta ]}^{p,n},\;Z_{t}^{p,n}=z_{[t/\delta
]}^{p,n},\;A_{t}^{p,n}=\sum_{m=0}^{[t/\delta ]}a_{m}^{p,n},\
K_{t}^{p,n}=\sum_{m=0}^{[t/\delta ]}k_{m}^{p,n},
\]
or $(\bar{Y}_{t}^{p,n},\bar{Z}_{t}^{p,n},\overline{A}_{t}^{p,n},\bar{K}%
_{t}^{p,n})$ by the implicit--explicit penalization scheme,
\[
\bar{Y}_{t}^{p,n}=\bar{y}_{[t/\delta ]}^{p,n},\;\bar{Z}_{t}^{p,n}=\bar{z}%
_{[t/\delta ]}^{p,n},\;\overline{A}_{t}^{p,n}=\sum_{m=0}^{[t/\delta ]}%
\overline{a}_{m}^{p,n},\ \bar{K}_{t}^{p,n}=\sum_{m=0}^{[t/\delta ]}\overline{%
k}_{m}^{p,n}.
\]
Let us notice that the laws of the solutions $(Y^{p},Z^{p},A^{p},K^{p})$ and
$(Y^{p,n},Z^{p,n},A^{p,n},K^{p,n})$ or $(\bar{Y}^{p,n},\bar{Z}^{p,n},%
\overline{A}^{p,n},\bar{K}^{p,n})$ to penalized BSDE depend only on $(%
\mathbf{P}_{B},\Gamma ^{-1}(\mathbf{P}_{B}),g,\Psi _{1}^{-1}(\mathbf{P}%
_{B}),\Psi _{2}^{-1}(\mathbf{P}_{B}))$ and $(\mathbf{P}_{B^{n}},\Gamma ^{-1}(%
\mathbf{P}_{B^{n}}),g,\Psi _{1}^{-1}(\mathbf{P}_{B^{n}}),\Psi _{2}^{-1}(%
\mathbf{P}_{B^{n}}))$ where $\mathbf{P}_{B}$(resp. $\mathbf{P}_{B^n}$) is the probability introduced by $B$(resp. $B^n$),
and $f^{-1}(\mathbf{P}_{B})$ (resp. $f^{-1}(\mathbf{P%
}_{B^{n}})$) is the law of $f(B)$ (resp. $f(B^{n})$) for $f=\Gamma
$, $\Psi _{1}$, $\Psi _{2}$. So if we concern the convergence in
law, we can consider these equations on any probability space.

By Donsker's theorem and Skorokhod representation theorem, there
exists a probability space, such that $\sup_{0\leq t\leq T}\left|
B_{t}^{n}-B_{t}\right| \rightarrow 0$, as $n\rightarrow \infty $, in $%
\mathbf{L}^{2}(\mathcal{F}_{T})$, since $\varepsilon _{k}$ is in $\mathbf{L}%
^{2+\delta }$. So we will work on this space with respect to the
filtration generated by $B^{n}$ and $B$, trying to prove the
convergence of solutions. Thanks to the convergence of $B^{n}$,
$(L^{n},U^{n})$ also converges to $(L,U)$. Then we have the
following result, which is based on the convergence results of
numerical solutions for BSDE (cf. \cite{bdm}, \cite{BDM2}) and
penalization method for reflected BSDE (Theorem \ref{p2b}).

\begin{proposition}
\label{conl2-pimp}Assume \ref{ass-ter} holds, the sequence
$(Y_{t}^{p,n},Z_{t}^{p,n})$ converges to $(Y_{t},Z_{t})$ in
following sense
\begin{equation}
\lim_{p\rightarrow \infty}\lim_{n\rightarrow \infty}E[\sup_{0\leq
t\leq T}\left| Y_{t}^{p,n}-Y_{t}\right| ^{2}+\int_{0}^{T}\left|
Z_{s}^{p,n}-Z_{s}\right| ^{2}ds]\rightarrow 0, \label{convl2-np}
\end{equation}
and for $0\leq t\leq T$, $A_{t}^{p,n}-K_{t}^{p,n}\rightarrow A_{t}-K_{t}$ in $\mathbf{L}^{2}(%
\mathcal{F}_{t})$, as$\ n\rightarrow \infty $, $%
p\rightarrow \infty $.
\end{proposition}

\textbf{Proof. }Notice
\begin{eqnarray*}
E[\sup_{0\leq t\leq T}\left| Y_{t}^{p,n}-Y_{t}\right|
^{2}+\int_{0}^{T}\left| Z_{s}^{p,n}-Z_{s}\right| ^{2}ds] &\leq
&2E[\sup_{0\leq t\leq T}\left| Y_{t}^{p,n}-Y_{t}^{p}\right|
^{2}+\int_{0}^{T}\left| Z_{s}^{p,n}-Z_{s}^{p}\right| ^{2}ds] \\
&&+2E[\sup_{0\leq t\leq T}\left| Y_{t}^{p}-Y_{t}\right|
^{2}+\int_{0}^{T}\left| Z_{s}^{p}-Z_{s}\right| ^{2}ds].
\end{eqnarray*}
By the convergence results of numerical solutions for BSDE (cf.
\cite{bdm}, \cite{BDM2}), the first part tends to 0. For the second
part, it is a direct application of Theorem \ref{p2b} of the
penalization method. So we get (\ref{convl2-np}). For the increasing
processes, we have
\begin{eqnarray*}
E[((A_{t}^{p,n}-K_{t}^{p,n})-(A_{t}-K_{t}))^{2}] &\leq
&2E[((A_{t}^{p,n}-K_{t}^{p,n})-(A_{t}^{p}-K_{t}^{p}))^{2}] \\
&&+2E[((A_{t}^{p}-K_{t}^{p})-(A_{t}-K_{t}))^{2}] \\
&\leq &2E[((A_{t}^{p,n}-K_{t}^{p,n})-(A_{t}^{p}-K_{t}^{p}))^{2}]+\frac{C}{%
\sqrt{p}},
\end{eqnarray*}
in view of (\ref{conl2-p}). While for fixed $p$,
\begin{eqnarray*}
A_{t}^{p,n}-K_{t}^{p,n}
&=&Y_{0}^{p,n}-Y_{t}^{p,n}-\int_{0}^{t}g(s,Y_{s}^{p,n},Z_{s}^{p,n})ds+%
\int_{0}^{t}Z_{s}^{p,n}dB_{s}^{n}, \\
A_{t}^{p}-K_{t}^{p}
&=&Y_{0}^{p}-Y_{t}^{p}-\int_{0}^{t}g(s,Y_{s}^{p},Z_{s}^{p})ds+%
\int_{0}^{t}Z_{s}^{p}dB_{s},
\end{eqnarray*}
from Corollary 14 in \cite{BDM2}, we know that $\int_{0}^{\cdot
}Z_{s}^{p,n}dB_{s}^{n}$ converges to $\int_{0}^{\cdot }Z_{s}^{p}dB_{s}$ in $%
\mathbf{S}^{2}(0,T)$, as $n\rightarrow \infty $, then with the
Lipschitz
condition of $g$ and the convergence of $Y^{p,n}$, we get $%
(A_{t}^{p,n}-K_{t}^{p,n})\rightarrow (A_{t}-K_{t})$ in $\mathbf{L}^{2}(%
\mathcal{F}_{t})$, as $n\rightarrow \infty $, $p\rightarrow \infty $. $%
\square $

Now we consider the implicit--explicit penalization scheme. From
Proposition 5 in \cite{PX08}, we know that for implicit--explicit
scheme, the difference between this solution and the totally
implicit one depends on $\mu+p$ for fixed $p\in\mathbb{N}$. So we
have

\begin{proposition}
For any $p\in \mathbb{N}$, when $n\rightarrow \infty $,
\[
E[\sup_{0\leq t\leq T}\left|
\overline{Y}_{t}^{p,n}-Y_{t}^{p,n}\right| ^{2}+\int_{0}^{T}\left|
\overline{Z}_{s}^{p,n}-Z_{s}^{p,n}\right| ^{2}ds\rightarrow
0\mbox{,}
\]
with
$(\overline{A}_{t}^{p,n}-\overline{K}_{t}^{p,n})-(A_{t}^{p,n}-K_{t}^{p,n})\rightarrow
0 $ in $\mathbf{L}^{2}(\mathcal{F}_{t})$, for $0\leq t\leq T$.
\end{proposition}

\textbf{Proof.} The convergence of $(\overline{Y}_{t}^{p,n},\overline{Z}%
_{t}^{p,n})$ to $(Y_{t}^{p,n},Z_{t}^{p,n})$ is a direct consequence
of Proposition 5 in \cite{PX08}. More precisely, there exists a
constant $C$ which depends only on $\mu $ and $T$, such that
\[
E[\sup_{0\leq t\leq T}\left|
\overline{Y}_{t}^{p,n}-Y_{t}^{p,n}\right| ^{2}]+E\int_{0}^{T}\left|
\overline{Z}_{s}^{p,n}-Z_{s}^{p,n}\right| ^{2}ds\leq C\delta ^{2}.
\]
Then we consider the convergence of the increasing processes, notice
that for $0\leq t\leq T$,
\[
\overline{A}_{t}^{p,n}-\overline{K}_{t}^{p,n}=\overline{Y}_{0}^{p,n}-%
\overline{Y}_{t}^{p,n}-\int_{0}^{t}g(s,\overline{Y}_{s}^{p,n},\overline{Z}%
_{s}^{p,n})ds+\int_{0}^{t}\overline{Z}_{s}^{p,n}dB_{s}^{n},
\]
compare with $A_{t}^{p,n}-K_{t}^{p,n}=Y_{0}^{p,n}-Y_{t}^{p,n}-%
\int_{0}^{t}g(s,Y_{s}^{p,n},Z_{s}^{p,n})ds+\int_{0}^{t}Z_{s}^{p,n}dB_{s}^{n}$%
, thanks to the Lipschitz condition of $g$ and the convergence of $(\overline{Y}%
^{p,n},\overline{Z}^{p,n})$, we get $\overline{A}_{t}^{p,n}-\overline{K}%
_{t}^{p,n}\rightarrow A_{t}^{p,n}-K_{t}^{p,n}$, in $\mathbf{L}^{2}(\mathcal{F%
}_{t})$, as $n\rightarrow \infty $, for fixed $p $. So the result
follows. $\square $

\begin{remark}
From this proposition and Proposition \ref{conl2-pimp}, we get the
convergence of the implicit--explicit penalization scheme.
\end{remark}

Before going further, we prove an \textit{a-priori} estimation of
$(y^{p,n},z^{p,n},a^{p,n},k^{p,n})$. This result will help us to get
the convergence of reflected scheme, which will be discussed in the
next section. Throughout this paper, we use $C_{\phi ,\psi, \cdots
}$ to denote a constant which depends on $\phi $, $\psi $, $\cdots$.
Here $\phi $, $\psi $, $\cdots$ can be random variables or
stochastic processes.

\begin{lemma}
\label{est-imp2b1}For each $p\in \mathbb{N}$ and $\delta $ such that
$\delta (1+2\mu +2\mu ^{2})<1$, there exists a constant $c$ such
that
\[
E[\sup_{j}\left| y_{j}^{p,n}\right| ^{2}+\sum_{j=0}^{n}\left|
z_{j}^{p,n}\right| ^{2}\delta +\frac{1}{p\delta
}\sum_{j=0}^{n}\left| a_{j}^{p,n}\right| ^{2}+\frac{1}{p\delta
}\sum_{j=0}^{n}\left| k_{j}^{p,n}\right| ^{2}]\leq c\;C_{\xi
^{n},g,L^{n},U^{n}}.
\]
Here $C_{\xi ^{n},g,L^{n},U^{n}}$ depends on $\xi ^{n}$, $g(t,0,0)$, $%
(L^{n})^{+}$ and $(U^{n})^{-}$, while $c$ depends only on $\mu $ and $T$%
.
\end{lemma}

\textbf{Proof.} Recall (\ref{dis-p2b}), we apply 'discrete It\^{o} formula' (cf. \cite{MPX}) for $%
(y_{j}^{p,n})^{2}$, we get
\begin{eqnarray*}
E[\left| y_{j}^{p,n}\right| ^{2}+\sum_{i=j}^{n-1}\left| z_{i}^{p,n}\right|
^{2}\delta ] &\leq &E[\left| \xi ^{n}\right|
^{2}+2[\sum_{i=j}^{n-1}y_{i}^{p,n}\left|
g(t_{i},y_{i}^{p,n},z_{i}^{p,n})\right|\delta ] \\
&&+2E\sum_{i=j}^{n-1}(y_{i}^{p,n}\cdot a_{i}^{p,n}-y_{i}^{p,n}\cdot
k_{i}^{p,n}).
\end{eqnarray*}
Since $y_{i}^{p,n}\cdot a_{i}^{p,n}=-p\delta
((y_{i}^{p,n}-L_{i}^{n})^{-})^{2}+p\delta
L_{i}^{n}(y_{i}^{p,n}-L_{i}^{n})^{-}=\frac{1}{p\delta }%
a_{i}^{p,n}+L_{i}^{n}a_{i}^{p,n}$ and $y_{i}^{p,n}\cdot k_{i}^{p,n}=p\delta
((y_{i}^{p,n}-U_{i}^{n})^{+})^{2}+U_{i}^{n}p\delta
(y_{i}^{p,n}-U_{i}^{n})^{+}=\frac{1}{p\delta }%
k_{i}^{p,n}+U_{i}^{n}k_{i}^{p,n}$, we have
\begin{eqnarray*}
&&E[\left| y_{j}^{p,n}\right| ^{2}+\frac{1}{2}\sum_{i=j}^{n-1}\left|
z_{i}^{p,n}\right| ^{2}\delta ]+2E[\frac{1}{p\delta }%
\sum_{i=j}^{n-1}(a_{i}^{p,n})^{2}+\frac{1}{p\delta }%
\sum_{i=j}^{n-1}(k_{i}^{p,n})^{2}] \\
&\leq &E[\left| \xi ^{n}\right| ^{2}+\sum_{i=j}^{n-1}\left|
g(t_{i},0,0)\right| ^{2}\delta +(1+2\mu +2\mu ^{2})\sum_{i=j}^{n-1}\left|
y_{i}^{p,n}\right| ^{2}\delta
+2\sum_{i=j}^{n-1}(L_{i}^{n})^{+}a_{i}^{p,n}+2%
\sum_{i=j}^{n-1}(U_{i}^{n})^{-}k_{i}^{p,n}] \\
&\leq &E[\left| \xi ^{n}\right| ^{2}+\sum_{i=j}^{n-1}\left|
g(t_{i},0,0)\right| ^{2}\delta ]+(1+2\mu +2\mu ^{2})\delta
E\sum_{i=j}^{n-1}\left| y_{i}^{p,n}\right| ^{2}+\frac{1}{\alpha }%
E(\sum_{i=j}^{n-1}a_{i}^{p,n})^{2} \\
&&+\alpha E\sup_{j\leq i\leq n-1}((L_{i}^{n})^{+})^{2}+\frac{1}{\beta }%
E(\sum_{i=j}^{n-1}k_{i}^{p,n})^{2}+\beta E\sup_{j\leq i\leq
n-1}((U_{i}^{n})^{+})^{2}.
\end{eqnarray*}
Since $L^n$ and $U^n$ are approximations of It\^o processes, we can find a process $%
X_{j}^{n}$ of the form $X_{j}^{n}=X_{0}-\sum_{i=0}^{j-1}\sigma
_{i}\varepsilon _{i+1}^{n}\sqrt{\delta }+V_{j}^{+n}-V_{j}^{-n}$, where $%
V_{j}^{\pm n}$ are $\mathcal{G}_{j}^{n}$-adapted increasing processes with $%
E[|V_{n}^{+n}|^{2}+|V_{n}^{-n}|^{2}]<+\infty $, and $L_{j}^{n}\leq
X_{j}^{n}\leq U_{j}^{n}$ holds. Then we apply similar techniques of
stopping times as in the proof of Lemma 2 in \cite{LSM} for the
discrete case with $L_{j}^{n}\leq X_{j}^{n}\leq U_{j}^{n}$, we can
prove
\[
E(\sum_{i=j}^{n-1}a_{i}^{p,n})^{2}+E(\sum_{i=j}^{n-1}k_{i}^{p,n})^{2}\leq
3\mu (C_{\xi ^{n},g,X^{n}}+E\sum_{i=j}^{n-1}[\left| y_{i}^{p,n}\right|
^{2}+\left| z_{i}^{p,n}\right| ^{2}]\delta ).
\]
While $X^n$ can be controlled by $L^n$ and $U^n$, we can replace it
by $L^n$ and $U^n$. Set $\alpha =\beta =12\mu $ in the previous
inequality, with Lemma \ref {dis-gro}, we get
\[
\sup_{j}E[\left| y_{j}^{p,n}\right| ^{2}]+E[\sum_{i=0}^{n-1}\left|
z_{i}^{p,n}\right| ^{2}\delta ]+\frac{1}{p\delta }%
\sum_{i=0}^{n-1}(a_{i}^{p,n})^{2}+\frac{1}{p\delta }%
\sum_{i=0}^{n-1}(k_{i}^{p,n})^{2}]\leq cC_{\xi ^{n},g,L^{n},U^{n}}.
\]
We reconsider It\^o formula for $|y_{j}^{p,n}|^2$, the take $\sup_j$
before expectation. Using Burkholder-Davis-Gundy inequality for
martingale part $\sum_{i=0}^{j}
y_{j}^{p,n}z_{j}^{p,n}\sqrt{\delta}\varepsilon_{j+1}^n$, with
similar techniques, we get
\[
E[\sup_{j}\left| y_{j}^{p,n}\right| ^{2}]\leq C_{\xi
^{n},g,L^{n},U^{n}}+C_{\mu}E[\sum_{i=0}^{n-1}\left|
y_{i}^{p,n}\right| ^{2}\delta ]\leq C_{\xi
^{n},g,L^{n},U^{n}}+C_{\mu}T \sup_j E[| y_{j}^{p,n}| ^{2} ].
\]
It follows the desired results. $\square $

\section{Reflected Algorithms and their convergence}

\subsection{Reflected Schemes}

This type of numerical schemes is based on reflecting the solution
$y^{n}$ between two barriers by $a^{n}$ and $k^{n}$ directly. In
such a way the
discrete solution $y^{n}$ really stays between two barriers $L^{n}$ and $%
U^{n}$. After discretization of time interval, our discrete reflected BSDE
with two barriers on small interval $[t_{j},t_{j+1}]$, for $0\leq j\leq n-1$%
, is
\begin{equation}
y_{j}^{n}=y_{j+1}^{n}+g(t_{j},y_{j}^{n},z_{j}^{n})\delta
+a_{j}^{n}-k_{j}^{n}-z_{j}^{n}\sqrt{\delta }\varepsilon _{j+1}^{n},
\label{dis-rbsde2b}
\end{equation}
with terminal condition $y_{n}^{n}=\xi ^{n}$, and constraint and discrete
integral conditions hold:
\begin{eqnarray}
a_{j}^{n} &\geq &0,\;\;k_{j}^{n}\geq 0,\ \ a_{j}^{n}\cdot k_{j}^{n}=0,
\label{dis-rbsde2b-r} \\
\;\;L_{j}^{n} &\leq &y_{j}^{n}\leq
U_{j}^{n},\;%
\;(y_{j}^{n}-L_{j}^{n})a_{j}^{n}=(y_{j}^{n}-U_{j}^{n})k_{j}^{n}=0.  \nonumber
\end{eqnarray}
Note that, all terms in (\ref{dis-rbsde2b}) are $\mathcal{G}_{j}^{n}$%
-measurable except $y_{j+1}^{n}$ and $\varepsilon _{j+1}^{n}$.

The key point of our numerical schemes is how to solve $%
(y_{j}^{n},z_{j}^{n},a_{j}^{n},k_{j}^{n})$ from (\ref{dis-rbsde2b})
using the $\mathcal{G}_{j+1}^{n}$-measurable random variable
$y_{j+1}^{n}$ obtained in the preceding step. First $z_{j}^{n}$ is
obtained by
\[
z_{j}^{n}=E[y_{j+1}^{n}\varepsilon _{j+1}^{n}|\mathcal{G}_{j}^{n}]=\frac{1}{2%
\sqrt{\delta }}(y_{j+1}^{n}|_{\varepsilon
_{j}^{n}=1}-y_{j+1}^{n}|_{\varepsilon _{j}^{n}=-1}).
\]
Then (\ref{dis-rbsde2b}) with (\ref{dis-rbsde2b-r}) becomes
\begin{eqnarray}
y_{j}^{n} &=&E[y_{j+1}^{n}|\mathcal{G}_{j}^{n}]+g(t_{j},y_{j}^{n},z_{j}^{n})%
\delta +a_{j}^{n}-k_{j}^{n},\;a_{j}^{n}\geq 0,\;k_{j}^{n}\geq 0,
\label{dis-rbsde2b1} \\
L_{j}^{n} &\leq &y_{j}^{n}\leq
U_{j}^{n},\;%
\;(y_{j}^{n}-L_{j}^{n})a_{j}^{n}=(y_{j}^{n}-U_{j}^{n})k_{j}^{n}=0.\;
\nonumber
\end{eqnarray}
Set $\Theta (y):=y-g(t_{j},y,z_{j}^{n})\delta $. In view of
$\left\langle \Theta (y)-\Theta (y^{\prime }),y-y^{\prime
}\right\rangle \geq (1-\delta \mu )\left| y-y^{\prime }\right|
^{2}>0$, for $\delta $ small enough, we get that in such case
$\Theta (y)$ is strictly increasing in $y$. So

\begin{eqnarray*}
y &\geq &L_{j}^{n}\Longleftrightarrow \Theta (y)\geq \Theta (L_{j}^{n}), \\
y &\leq &U_{j}^{n}\Longleftrightarrow \Theta (y)\leq \Theta (U_{j}^{n}).
\end{eqnarray*}
Then implicit reflected scheme gives the results with $E[y_{j+1}^{n}|%
\mathcal{G}_{j}^{n}]=\frac{1}{2}(y_{j+1}^{n}|_{\varepsilon
_{j}^{n}=1}+y_{j+1}^{n}|_{\varepsilon _{j}^{n}=-1})$ as follows
\begin{eqnarray*}
y_{j}^{n} &=&\Theta ^{-1}(E[y_{j+1}^{n}|\mathcal{G}%
_{j}^{n}]+a_{j}^{n}-k_{j}^{n}), \\
a_{j}^{n} &=&\left( E[y_{j+1}^{n}|\mathcal{G}%
_{j}^{n}]+g(t_{j},L_{j}^{n},z_{j}^{n})\delta -L_{j}^{n}\right) ^{-}, \\
k_{j}^{n} &=&\left( E[y_{j+1}^{n}|\mathcal{G}%
_{j}^{n}]+g(t_{j},U_{j}^{n},z_{j}^{n})\delta -U_{j}^{n}\right) ^{+},
\end{eqnarray*}
on the set $\{L_{j}^{n}<U_{j}^{n}\}$, then we know that $\{y_{j}^{n}-L_{j}^{n}=0%
\}$ and $\{y_{j}^{n}-U_{j}^{n}=0\}$ are disjoint. So with $%
\;(y_{j}^{n}-L_{j}^{n})a_{j}^{n}=(y_{j}^{n}-U_{j}^{n})k_{j}^{n}=0$, we have $%
a_{j}^{n}\cdot k_{j}^{n}=0$. On the set $\{L_{j}^{n}=U_{j}^{n}\}$,
 we get $a^n_j=(I_j^n)^+$ and $k^n_j=(I_j^n)^-$ by defintion, where $I^n_j :=E[y_{j+1}^{n}|\mathcal{G}%
_{j}^{n}]+g(t_{j},L_{j}^{n},z_{j}^{n})\delta -L_{j}^{n}$. So
automatically $a_{j}^{n}\cdot k_{j}^{n}=0$.

Our explicit reflected scheme is introduced by replacing $y_{j}^{n}$
in $g$ by $E[\bar{y}_{j+1}^{n}|\mathcal{G}_{j}^{n}]$ in
(\ref{dis-rbsde2b1}). So we get the following equation,
\begin{eqnarray}
\bar{y}_{j}^{n} &=&E[\bar{y}_{j+1}^{n}|\mathcal{G}_{j}^{n}]+g(t_{j},E[\bar{y}%
_{j+1}^{n}|\mathcal{G}_{j}^{n}],\bar{z}_{j}^{n})\delta +\overline{a}_{j}^{n}-%
\overline{k}_{j}^{n},\;\overline{a}_{j}^{n}\geq 0,\;\overline{k}_{j}^{n}\geq
0,  \label{pena-exp1} \\
L_{j}^{n} &\leq &\bar{y}_{j}^{n}\leq U_{j}^{n},\;\;(\bar{y}%
_{j}^{n}-L_{j}^{n})\overline{a}_{j}^{n}=(\bar{y}_{j}^{n}-U_{j}^{n})\overline{%
k}_{j}^{n}=0.\;  \nonumber
\end{eqnarray}
Then with $E[\overline{y}_{j+1}^{n}|\mathcal{G}_{j}^{n}]=\frac{1}{2}(%
\overline{y}_{j+1}^{n}|_{\varepsilon _{j}^{n}=1}+\overline{y}%
_{j+1}^{n}|_{\varepsilon _{j}^{n}=-1})$, we get the solution
\begin{eqnarray}
\overline{y}_{j}^{n} &=&E[\bar{y}_{j+1}^{n}|\mathcal{G}_{j}^{n}]+g(t_{j},E[%
\bar{y}_{j+1}^{n}|\mathcal{G}_{j}^{n}],\overline{z}_{j}^{n})\delta +%
\overline{a}_{j}^{n}-\overline{k}_{j}^{n},  \nonumber \\
\overline{a}_{j}^{n} &=&\left( E[\bar{y}_{j+1}^{n}|\mathcal{G}%
_{j}^{n}]+g(t_{j},E[\bar{y}_{j+1}^{n}|\mathcal{G}_{j}^{n}],\bar{z}%
_{j}^{n})\delta -L_{j}^{n}\right) ^{-},  \label{p-impli} \\
\overline{k}_{j}^{n} &=&\left( E[\bar{y}_{j+1}^{n}|\mathcal{G}%
_{j}^{n}]+g(t_{j},E[\bar{y}_{j+1}^{n}|\mathcal{G}_{j}^{n}],\bar{z}%
_{j}^{n})\delta -U_{j}^{n}\right) ^{+}.  \nonumber
\end{eqnarray}

\subsection{Convergence of Reflected Implicit Schemes}

Now we study the convergence of Reflected Schemes. For implicit
reflected scheme, we denote
\[
Y_{t}^{n}=y_{[t/\delta ]}^{n},\;Z_{t}^{n}=z_{[t/\delta
]}^{n},\;\;A_{t}^{n}=\sum_{i=0}^{[t/\delta
]}a_{i}^{n},\;\;K_{t}^{n}=\sum_{i=0}^{[t/\delta ]}k_{i}^{n},
\]
and for explicit reflected scheme
\[
\bar{Y}_{t}^{n}=\bar{y}_{[t/\delta ]}^{n},\;\bar{Z}_{t}^{n}=\bar{z}%
_{[t/\delta ]}^{n},\;\;\overline{A}_{t}^{n}=\sum_{i=0}^{[t/\delta ]}%
\overline{a}_{i}^{n},\;\;\bar{K}_{t}^{n}=\sum_{i=0}^{[t/\delta ]}\overline{k}%
_{i}^{n}.
\]

First we prove an estimation result for $(y^{n},z^{n},a^{n},k^{n})$.

\begin{lemma}
\label{est-imp2b2}For $\delta $ such that $\delta (1+2\mu +2\mu
^{2})<1$, there exists a constant $c$ depending only on $\mu $ and
$T$ such that
\[
E[\sup_{j}\left| y_{j}^{n}\right| ^{2}+\sum_{j=0}^{n-1}\left|
z_{j}^{n}\right| ^{2}\delta +\left| \sum_{j=0}^{n-1}a_{j}^{n}\right|
^{2}+\left| \sum_{j=0}^{n-1}k_{j}^{n}\right| ^{2}]\leq c\;C_{\xi
^{n},g,L^{n},U^{n}}.
\]
\end{lemma}

\textbf{Proof.} First we consider the estimation of $a_{i}^{n}$ and $%
k_{i}^{n}$. In view of $L_{j}^{n}\leq Y_{j}^{n}\leq U_{j}^{n}$, we have that
\begin{eqnarray}
a_{j}^{n} &\leq &\left( E[L_{j+1}^{n}|\mathcal{G}%
_{j}^{n}]+g(t_{j},L_{j}^{n},z_{j}^{n})\delta -L_{j}^{n}\right) ^{-}=\delta
(l_{j}+g(t_{j},L_{j}^{n},z_{j}^{n}))^{-},  \label{increase} \\
k_{j}^{n} &\leq &\left( E[U_{j+1}^{n}|\mathcal{G}
_{j}^{n}]+g(t_{j},U_{j}^{n},z_{j}^{n})\delta -U_{j}^{n}\right)
^{+}=\delta \left( u_{j}+g(t_{j},U_{j}^{n},z_{j}^{n})\right) ^{+}.
\nonumber
\end{eqnarray}
We consider following discrete BSDEs with $\widehat{y}_{n}^{n}=\widetilde{y}%
_{n}^{n}=\xi ^{n}$,
\begin{eqnarray*}
\widehat{y}_{j}^{n} &=&\widehat{y}_{j+1}^{n}+[g(t_{j},\widehat{y}_{j}^{n},%
\widehat{z}_{j}^{n})+(l_{j})^{-}+g(t_{j},L_{j}^{n},\widehat{z}%
_{j}^{n})^{-}]\delta -\widehat{z}_{j}^{n}\sqrt{\delta }\varepsilon
_{j+1}^{n}, \\
\widetilde{y}_{j}^{n} &=&\widetilde{y}_{j+1}^{n}+[g(t_{j},\widetilde{y}%
_{j}^{n},\widetilde{z}_{j}^{n})-(u_{j})^{+}-g(t_{j},U_{j}^{n},\widetilde{z}%
_{j}^{n})^{+}]\delta -\widetilde{z}_{j}^{n}\sqrt{\delta }\varepsilon
_{j+1}^{n}.
\end{eqnarray*}
Thanks to discrete comparison theorem in \cite{MPX}, we have $\widetilde{y}%
_{j}^{n}\leq y_{j}^{n}\leq \widehat{y}_{j}^{n}$, so
\begin{equation}
E[\sup_{j}\left| y_{j}^{n}\right| ^{2}]\leq \max \{E[\sup_{j}\left|
\widetilde{y}_{j}^{n}\right| ^{2}],E[\sup_{j}\left| \widehat{y}%
_{j}^{n}\right| ^{2}]\}\leq c\;C_{\xi ^{n},g,L^{n},U^{n}}.
\label{est-y-d2b}
\end{equation}
The last inequality follows from estimations of discrete solution of
classical BSDE $(\widehat{y}_{j}^{n})^{2}$ and
$(\widetilde{y}_{j}^{n})^{2}$, which is obtained by It\^o formulae
and the discrete
Gronwall inequality in Lemma \ref{dis-gro}. For $z^n_j$, we use 'discrete It\^{o} formula' (cf. \cite{MPX}) again for $%
(y_{j}^{n})^{2}$, and get
\begin{eqnarray*}
E\left| y_{j}^{n}\right| ^{2}+\sum_{i=j}^{n-1}\left| z_{i}^{n}\right|
^{2}\delta &=&E[\left| \xi ^{n}\right|
^{2}+2\sum_{i=j}^{n-1}y_{i}^{n}g(t_{i},y_{i}^{n},z_{i}^{n})\delta
+2\sum_{i=j}^{n-1}y_{i}^{n}a_{i}^{n}-2\sum_{i=j}^{n-1}y_{i}^{n}k_{i}^{n}] \\
&\leq &E[\left| \xi ^{n}\right| ^{2}+\sum_{i=j}^{n-1}\left|
g(t_{i},0,0)\right| ^{2}\delta +\delta (1+2\mu +2\mu
^{2})\sum_{i=j}^{n-1}\left| y_{i}^{n}\right| ^{2}+\frac{1}{2}%
\sum_{i=j}^{n-1}\left| z_{i}^{n}\right| ^{2}\delta ] \\
&&+\alpha E[\sup_{j}((L_{j}^{n})^{+})^{2}+\sup_{j}((U_{j}^{n})^{-})^{2}]+%
\frac{1}{\alpha }E[(\sum_{j=i}^{n-1}a_{j}^{n})^{2}+(%
\sum_{j=i}^{n-1}k_{j}^{n})^{2}],
\end{eqnarray*}
using $(y_{i}^{n}-L_{i}^{n})a_{i}^{n}=0$ and $%
(y_{i}^{n}-U_{i}^{n})k_{i}^{n}=0$. And from (\ref{increase}), we
have
\begin{eqnarray*}
E(\sum_{j=i}^{n-1}a_{j}^{n})^{2} &\leq &4\delta
E\sum_{j=i}^{n-1}[(l_{j})^{2}+g(t_{i},0,0)^{2}+\mu \left| L_{j}^{n}\right|
^{2}+\mu \left| z_{j}^{n}\right| ^{2}], \\
E(\sum_{j=i}^{n-1}k_{j}^{n})^{2} &\leq &4\delta
E\sum_{j=i}^{n-1}[(u_{j})^{2}+g(t_{i},0,0)^{2}+\mu \left| U_{j}^{n}\right|
^{2}+\mu \left| z_{j}^{n}\right| ^{2}].
\end{eqnarray*}
Set $\alpha =32\mu $, it follows
\begin{eqnarray*}
E[\left| y_{j}^{n}\right| ^{2}+\frac{1}{4}\sum_{i=j}^{n-1}\left|
z_{i}^{n}\right| ^{2}\delta ] &\leq &E[\left| \xi ^{n}\right| ^{2}+(1+\frac{1%
}{8\mu ^{2}})\sum_{i=j}^{n-1}\left| g(t_{i},0,0)\right| ^{2}\delta ]+\delta
(1+2\mu +2\mu ^{2})\sum_{i=j}^{n-1}\left| y_{i}^{n}\right| ^{2} \\
&&+32\mu ^{2}E[\sup_{j}((L_{j}^{n})^{+})^{2}+\sup_{j}((U_{j}^{n})^{-})^{2}]+%
\frac{1}{8\mu ^{2}}E\sum_{j=i}^{n-1}[(l_{j})^{2}+(u_{j})^{2}] \\
&&+\frac{1}{8}\delta E\sum_{j=i}^{n-1}[\left| L_{j}^{n}\right| ^{2}+\left|
U_{j}^{n}\right| ^{2}]
\end{eqnarray*}
With (\ref{est-y-d2b}), we obtain $\sum_{i=0}^{n-1}\left|
z_{i}^{n}\right| ^{2}\delta \leq cC_{\xi ^{n},g,L^{n},U^{n}}$. Then
applying these estimations to (\ref{increase}), we obtain desired
results. $\square $

With arguments similar to those precede Proposition
\ref{conl2-pimp}, the laws of the solutions $(Y,Z,A,K)$ and
$(Y^{n},Z^{n},A^{n},K^{n})$ or $(\bar{Y}^{n},\bar{Z}^{n},%
\overline{A}^{n},\bar{K}^{n})$ to reflected BSDE depend only on $(%
\mathbf{P}_{B},\Gamma ^{-1}(\mathbf{P}_{B}),g$, $\Psi _{1}^{-1}(\mathbf{P}%
_{B}),\Psi _{2}^{-1}(\mathbf{P}_{B}))$ and $(\mathbf{P}_{B^{n}},\Gamma ^{-1}(%
\mathbf{P}_{B^{n}}),g,\Psi _{1}^{-1}(\mathbf{P}_{B^{n}}),\Psi _{2}^{-1}(%
\mathbf{P}_{B^{n}}))$ where $f^{-1}(\mathbf{P}_{B})$ (resp. $f^{-1}(\mathbf{P%
}_{B^{n}})$) is the law of $f(B)$ (resp. $f(B^{n})$) for $f=\Gamma
$, $\Psi _{1}$, $\Psi _{2}$. So if we concern the convergence in
law, we can consider these equations on any probability space.

From Donsker's theorem and Skorokhod representation theorem, there
exists a probability space satisfying $\sup_{0\leq t\leq T}\left|
B_{t}^{n}-B_{t}\right|
\rightarrow 0$, as $n\rightarrow \infty $, in $\mathbf{L}^{2}(\mathcal{F}%
_{T})$, since $\varepsilon _{k}$ is in $\mathbf{L}^{2+\delta }$. And
it is sufficient for us to prove convergence results in this
probability space. Our convergence result for the implicit reflected
scheme is as follows:

\begin{theorem}
\label{conv-r2b}Under Assumption \ref{ass-ter} and suppose moreover
that $g$ satisfies Lipschitz
condition (\ref{Lip}), we have when $n\rightarrow +\infty ,$%
\begin{equation}
E[\sup_{t}|Y_{t}^{n}-Y_{t}|^{2}]+E\int_{0}^{T}|Z_{t}^{n}-Z_{t}|^{2}dt%
\rightarrow 0,  \label{conv-exp-R2b1}
\end{equation}
and $A_{t}^{n}-K_{t}^{n}\rightarrow A_{t}-K_{t}$ in $\mathbf{L}^{2}(\mathcal{%
F}_{t})$, for $0\leq t\leq T$.
\end{theorem}

\textbf{Proof.} The proof is done in three steps.

In the first step, we consider the difference between discrete
solutions of reflect implicit scheme and of penalization implicit
scheme introduce in section 4.1 and section 3.1, respectively. More
precisely, we will prove that for each $p$,
\begin{equation}
E[\sup_{j}|y_{j}^{n}-y_{j}^{p,n}|^{2}]+\delta
E\sum_{j=0}^{n-1}|z_{j}^{n}-z_{j}^{p,n}|^{2}\leq cC_{\xi ^{n},g,L^{n},U^{n}}%
\frac{1}{\sqrt{p}}.  \label{diff-2bpr1}
\end{equation}
Here $c$ only depends on $\mu $ and $T$. From (\ref{dis-p2b}) and (\ref
{dis-rbsde2b}), applying 'discrete It\^{o} formula' (cf. \cite{MPX}) to $%
(y_{j}^{n}-y_{j}^{p,n})^{2}$, we get
\begin{eqnarray*}
&&E\left| y_{j}^{n}-y_{j}^{p,n}\right| ^{2}+\delta
E\sum_{i=j}^{n-1}|z_{i}^{n}-z_{i}^{p,n}|^{2} \\
&=&2E%
\sum_{i=j}^{n-1}[(y_{i}^{n}-y_{i}^{p,n})(g(t_{i},y_{i}^{n},z_{i}^{n})-g(t_{i},y_{i}^{p,n},z_{i}^{p,n}))\delta ]
\\
&&+2E\sum_{i=j}^{n-1}[(y_{i}^{n}-y_{i}^{p,n})(a_{i}^{n}-a_{i}^{p,n})]-2E%
\sum_{i=j}^{n-1}[(y_{i}^{n}-y_{i}^{p,n})(k_{i}^{n}-k_{i}^{p,n})]
\end{eqnarray*}
From (\ref{dis-rbsde2b-r}), we have
\begin{eqnarray*}
(y_{i}^{n}-y_{i}^{p,n})(a_{i}^{n}-a_{i}^{p,n})& =&
(y_{i}^{n}-L_i^n)a_{i}^{n}-(y_{i}^{p,n}-L_i^n)a_{i}^{n}-(y_{i}^{n}-L_i^n)a_{i}^{p,n}+(y_{i}^{p,n}-L_i^n)a_{i}^{p,n}\\
&\leq &(y_{i}^{p,n}-L_{i}^{n})^{-}a_{i}^{n}-((y_{i}^{p,n}-L_i^n)^-)^2, \\
&\leq &(y_{i}^{p,n}-L_{i}^{n})^{-}a_{i}^{n},
\end{eqnarray*}
Similarly we have $(y_{i}^{n}-y_{i}^{p,n})(k_{i}^{n}-k_{i}^{p,n})
\geq -(y_{i}^{p,n}-U_{i}^{n})k_{i}^{n}$. By (\ref{increase}) and the
Lipschitz property of $g$, it follows
\begin{eqnarray*}
&&E\left| y_{j}^{n}-y_{j}^{p,n}\right| ^{2}+\frac{\delta }{2}%
E\sum_{i=j}^{n-1}|z_{i}^{n}-z_{i}^{p,n}|^{2} \\
&\leq &(2\mu +2\mu ^{2})\delta
E\sum_{i=j}^{n-1}[(y_{i}^{n}-y_{i}^{p,n})^{2}]+2E%
\sum_{i=j}^{n-1}[(y_{i}^{p,n}-L_{i}^{n})^{-}a_{i}^{n}+(y_{i}^{p,n}-U_{i}^{n})^{+}k_{i}^{n}]
\\
&\leq &(2\mu +2\mu ^{2})\delta
E\sum_{i=j}^{n-1}[(y_{i}^{n}-y_{i}^{p,n})^{2}]+2\left( \delta
E\sum_{i=j}^{n-1}((y_{i}^{p,n}-L_{i}^{n})^{-})^{2}\right) ^{\frac{1}{2}%
}\left( \delta
E\sum_{i=j}^{n-1}((l_{j}+g(t_{i},L_{j}^{n},z_{j}^{n}))^{-})^{2}\right) ^{%
\frac{1}{2}} \\
&&+2\left( \delta E\sum_{i=j}^{n-1}((y_{i}^{p,n}-U_{i}^{n})^{+})^{2}\right)
^{\frac{1}{2}}\left( \delta
E\sum_{i=j}^{n-1}((u_{j}+g(t_{i},U_{j}^{n},z_{j}^{n}))^{+})^{2}\right) ^{\frac{1}{2}} \\
&=&(2\mu +2\mu ^{2})\delta E\sum_{i=j}^{n-1}[(y_{i}^{n}-y_{i}^{p,n})^{2}]+%
\frac{2}{\sqrt{p}}\left( \frac{1}{p\delta }%
E\sum_{i=j}^{n-1}(a_{i}^{p,n})^{2}\right) ^{\frac{1}{2}}\left( \delta
E\sum_{i=j}^{n-1}((l_{j}+g(t_{i},L_{j}^{n},z_{j}^{n}))^{-})^{2}\right) ^{%
\frac{1}{2}} \\
&&+\frac{2}{\sqrt{p}}\left( \frac{1}{p\delta }%
E\sum_{i=j}^{n-1}(k_{i}^{p,n})^{2}\right) ^{\frac{1}{2}}\left(
\delta
E\sum_{i=j}^{n-1}((u_{j}+g(t_{i},U_{j}^{n},z_{j}^{n}))^{+})^{2}\right)
^{\frac{1}{2}} .
\end{eqnarray*}
Then by estimation results in Lemma \ref{est-imp2b1}, Lemma
\ref{est-imp2b2} and discrete Gronwall inequality in Lemma
\ref{dis-gro}, we get
\[
\sup_{j}E\left| y_{j}^{n}-y_{j}^{p,n}\right| ^{2}+\delta
E\sum_{i=0}^{n-1}|z_{i}^{n}-z_{i}^{p,n}|^{2}\leq c\;C_{\xi ^{n},g,L^{n},U^{n}}%
\frac{1}{\sqrt{p}}.
\]
Apply B-D-G inequality, we obtain (\ref{diff-2bpr1}).

In the second step, we want to prove (\ref{conv-exp-R2b1}). We have
\begin{eqnarray*}
&&E[\sup_{t}|Y_{t}^{n}-Y_{t}|^{2}]+E[\int_{0}^{T}|Z_{t}^{n}-Z_{t}|^{2}dt] \\
&\leq
&3E[\sup_{t}|Y_{t}^{p}-Y_{t}|^{2}+\int_{0}^{T}|Z_{t}^{p}-Z_{t}|^{2}dt]+3E[%
\sup_{t}|Y_{t}^{n}-Y_{t}^{p,n}|^{2}+%
\int_{0}^{T}|Z_{t}^{n}-Z_{t}^{p,n}|^{2}dt] \\
&&+3E[\sup_{t}|Y_{t}^{p}-Y_{t}^{p,n}|^{2}+%
\int_{0}^{T}|Z_{t}^{p}-Z_{t}^{p,n}|^{2}dt] \\
&\leq &3Cp^{-\frac{1}{2}}+cC_{\xi ^{n},g,L^{n},U^{n}}p^{-\frac{1}{2}%
}+3E[\sup_{t}|Y_{t}^{p}-Y_{t}^{p,n}|^{2}+%
\int_{0}^{T}|Z_{t}^{p}-Z_{t}^{p,n}|^{2}dt],
\end{eqnarray*}
in view of (\ref{diff-2bpr1}) and Theorem \ref{p2b}. For fixed
$p>0$, by convergence results of numerical algorithms for BSDE,
(Theorem 12 in \cite {BDM2} and Theorem 2 in \cite{PX08}), we know
that the last two terms converge to $0$, as $\delta \rightarrow 0$.
And when $\delta $ is small enough, $C_{\xi ^{n},g,L^{n},U^{n}}$ is
dominated by $\xi ^{n}$, $g$, $L$ and $U$. This implies that we can
choose suitable $\delta $ such that the right hand side is as small
as we want, so (\ref{conv-exp-R2b1}) follows.

In the last step, we consider the convergence of $(A^{n},K^{n})$. Recall
that for $0\leq t\leq T$,
\begin{eqnarray*}
A_{t}^{n}-K_{t}^{n}
&=&Y_{0}^{n}-Y_{t}^{n}-\int_{0}^{t}g(s,Y_{s}^{n},Z_{s}^{n})ds+%
\int_{0}^{t}Z_{s}^{n}dB_{s}^{n}, \\
A_{t}^{p,n}-K_{t}^{p,n}
&=&Y_{0}^{p,n}-Y_{t}^{p,n}-\int_{0}^{t}g(s,Y_{s}^{p,n},Z_{s}^{p,n})ds+%
\int_{0}^{t}Z_{s}^{p,n}dB_{s}^{n}.
\end{eqnarray*}
By (\ref{diff-2bpr1}) and Lipschitz condition of $g$, we get
\[
E[\left| (A_{t}^{n}-K_{t}^{n})-(A_{t}^{p,n}-K_{t}^{p,n})\right|
^{2}]\leq c\;C_{\xi ^{n},g,L^{n},U^{n}}\frac{1}{\sqrt{p}}.
\]
Since
\begin{eqnarray*}
E[\left| (A_{t}^{n}-K_{t}^{n})-(A_{t}-K_{t})\right| ^{2}] &\leq &3E[\left|
(A_{t}^{n}-K_{t}^{n})-(A_{t}^{p,n}-K_{t}^{p,n})\right| ^{2}]+3E[\left|
(A_{t}^{p}-K_{t}^{p})-(A_{t}-K_{t})\right| ^{2}] \\
&&+3E[\left| (A_{t}^{p}-K_{t}^{p})-(A_{t}^{p,n}-K_{t}^{p,n})\right| ^{2}] \\
&\leq &c(C_{\xi ^{n},g,L^{n},U^{n}}+C_{\xi ,g,L,U})\frac{1}{\sqrt{p}}%
+3E[\left| (A_{t}^{p}-K_{t}^{p})-(A_{t}^{p,n}-K_{t}^{p,n})\right| ^{2}],
\end{eqnarray*}
with similar techniques, we obtain $E[\left|
(A_{t}^{n}-K_{t}^{n})-(A_{t}-K_{t})\right| ^{2}]\rightarrow 0$. Here the
fact that $(A_{t}^{p,n}-K_{t}^{p,n})$ converges to $(A_{t}^{p}-K_{t}^{p})$
for fixed $p$ follows from the convergence results of $%
(Y_{t}^{p,n},Z_{t}^{p,n})$ to $(Y_{t}^{p},Z_{t}^{p})$. $\square $

\subsection{Convergence of Reflected Explicit Scheme}
Then we study the convergence of explicit reflected scheme. Before going
further, we need an estimation result for $(\overline{y}^{n},%
\overline{z}^{n},\overline{a}^{n},\overline{k}^{n})$.

\begin{lemma}
\label{est-expRBSDE2d}For $\delta $ such that $\delta (\frac{9}{4}+2\mu
+4\mu ^{2})<1$, there exists a constant $c$ depending only on $\mu $ and $T$%
, such that
\[
E[\sup_{j}\left| \overline{y}_{j}^{n}\right|
^{2}]+E[\sum_{j=0}^{n-1}\left|
\overline{z}_{j}^{n}\right| ^{2}\delta +\left| \sum_{j=0}^{n-1}\overline{k}%
_{j}^{n}\right| ^{2}+\left|
\sum_{j=0}^{n-1}\overline{a}_{j}^{n}\right| ^{2}]\leq c\;C_{\xi
^{n},g,L^{n},U^{n}}.
\]
\end{lemma}

\textbf{Proof. }We recall the explicit reflected scheme, which is:
\begin{eqnarray}
\bar{y}_{j}^{n} &=&\bar{y}_{j+1}^{n}+g(t_{j},E[\bar{y}%
_{j+1}^{n}|\mathcal{G}_{j}^{n}],\bar{z}_{j}^{n})\delta +\overline{a}_{j}^{n}-%
\overline{k}_{j}^{n}-\bar{z}^n_j\sqrt{\delta}\varepsilon_{j+1}^n,\;\overline{a}_{j}^{n}\geq
0,\;\overline{k}_{j}^{n}\geq
0, \label{explicit}\\
L_{j}^{n} &\leq &\bar{y}_{j}^{n}\leq U_{j}^{n},\;\;(\bar{y}%
_{j}^{n}-L_{j}^{n})\overline{a}_{j}^{n}=(\bar{y}_{j}^{n}-U_{j}^{n})\overline{%
k}_{j}^{n}=0.\;\nonumber
\end{eqnarray}
Then we have
\begin{eqnarray}
|\bar{y}_{j}^{n}|^{2} &=&|\overline{y}_{j+1}^{n}|^{2}-|\bar{z}%
_{j}^{n}|^{2}\delta +2\bar{y}_{j+1}^{n}\cdot g(t_{j},E[\overline{y}%
_{j+1}^{n}|\mathcal{G}_{j}^{n}],\overline{z}_{j}^{n})\delta +2\bar{y}%
_{j}^{n}\cdot \overline{a}_{j}^{n}-2\bar{y}_{j}^{n}\cdot \overline{k}%
_{j}^{n} \label{sqr-exp}\\
&&+|g(t_{j},E[\overline{y}_{j+1}^{n}|\mathcal{G}_{j}^{n}],\overline{z}%
_{j}^{n})|^{2}\delta ^{2}-(\overline{a}_{j}^{n})^{2}-(\overline{k}%
_{j}^{n})^{2}-2\bar{y}_{j}^{n}\bar{z}^n_j\sqrt{\delta}\varepsilon_{j+1}^n\nonumber\\
&&+ 2g(t_{j},E[\overline{y}%
_{j+1}^{n}|\mathcal{G}_{j}^{n}],\overline{z}_{j}^{n})\bar{z}^n_j\delta\sqrt{\delta}\varepsilon_{j+1}^n
-2(\overline{a}_{j}^{n}-%
\overline{k}_{j}^{n})\bar{z}^n_j\sqrt{\delta}\varepsilon_{j+1}^n.\nonumber
\end{eqnarray}
In view of $(\overline{y}_{j}^{n}-L_{j}^{n})\overline{a}_{j}^{n}=(\bar{y}%
_{j}^{n}-U_{j}^{n})\overline{k}_{j}^{n}=0$, $\overline{a}_{j}^{n}$ and $%
\overline{k}_{j}^{n}\geq 0$, and taking expectation, we have
\begin{eqnarray*}
E|\bar{y}_{j}^{n}|^{2}+E|\bar{z}_{j}^{n}|^{2}\delta &\leq &E|\overline{y}%
_{j+1}^{n}|^{2}+2E[\bar{y}_{j+1}^{n}\cdot g(t_{j},E[\overline{y}_{j+1}^{n}|%
\mathcal{G}_{j}^{n}],\overline{z}_{j}^{n})]\delta +2E[(L_{j}^{n})^{+}\cdot
\overline{a}_{j}^{n}]+E[(U_{j}^{n})^{-}\cdot \overline{k}_{j}^{n}] \\
&&+E[|g(t_{j},E[\overline{y}_{j+1}^{n}|\mathcal{G}_{j}^{n}],\overline{z}%
_{j}^{n})|^{2}\delta ^{2}] \\
&\leq &E|\overline{y}_{j+1}^{n}|^{2}+(\delta +3\delta
^{2})E[|g(t_{j},0,0)|^{2}]+(\frac{1}{4}\delta +3\mu ^{2}\delta ^{2})E[(%
\overline{z}_{j}^{n})^{2}] \\
&&+\delta (1+2\mu +4\mu ^{2}+3\mu ^{2}\delta )E|\overline{y}%
_{j+1}^{n}|^{2}+2E[(L_{j}^{n})^{+}\cdot \overline{a}%
_{j}^{n}]+E[(U_{j}^{n})^{-}\cdot \overline{k}_{j}^{n}]
\end{eqnarray*}
Taking the sum for $j=i,\cdots ,n-1$ yields
\begin{eqnarray}
&&E|\bar{y}_{i}^{n}|^{2}+\frac{1}{2}\sum_{j=i}^{n-1}E|\bar{z}%
_{j}^{n}|^{2}\delta  \label{est-imp2b3} \\
&\leq &E|\xi ^{n}|^{2}+(\delta +3\delta
^{2})E\sum_{j=i}^{n-1}[|g(t_{j},0,0)|^{2}]+\delta (1+2\mu +4\mu ^{2}+3\mu
^{2}\delta )E\sum_{j=i}^{n-1}|\overline{y}_{j+1}^{n}|^{2}  \nonumber \\
&&+\alpha E[\sup_{j}((L_{j}^{n})^{+})^{2}+\sup_{j}((U_{j}^{n})^{+})^{2}]+%
\frac{1}{\alpha }E[(\sum_{j=i}^{n-1}\overline{a}_{j}^{n})^{2}+(%
\sum_{j=i}^{n-1}\overline{k}_{j}^{n})^{2}],  \nonumber
\end{eqnarray}
where $\alpha $ is a constant to be decided later. From (\ref{p-impli}), we
have
\begin{eqnarray*}
\overline{a}_{j}^{n} &\leq &\left( E[L_{j+1}^{n}|\mathcal{G}%
_{j}^{n}]+g(t_{j},E[\bar{y}_{j+1}^{n}|\mathcal{G}_{j}^{n}],\bar{z}%
_{j}^{n})\delta -L_{j}^{n}\right) ^{-}=(l_{j}+g(t_{j},E[\bar{y}_{j+1}^{n}|%
\mathcal{G}_{j}^{n}],\bar{z}_{j}^{n}))^{-}\delta , \\
\overline{k}_{j}^{n} &\leq &\left( E[U_{j+1}^{n}|\mathcal{G}%
_{j}^{n}]+g(t_{j},E[\bar{y}_{j+1}^{n}|\mathcal{G}_{j}^{n}],\bar{z}%
_{j}^{n})\delta -U_{j}^{n}\right) ^{+}=(u_{j}+g(t_{j},E[\bar{y}_{j+1}^{n}|%
\mathcal{G}_{j}^{n}],\bar{z}_{j}^{n}))^{-}\delta .
\end{eqnarray*}
Then we get
\begin{eqnarray}
E(\sum_{j=i}^{n-1}\overline{a}_{j}^{n})^{2} &\leq &4\delta
E\sum_{j=i}^{n-1}[(l_{j})^{2}+g(t_{j},0,0)^{2}+\mu ^{2}(E[\bar{y}_{j+1}^{n}|%
\mathcal{G}_{j}^{n}])^{2}+\mu ^{2}(\bar{z}_{j}^{n})^{2}],  \label{est-2b} \\
E(\sum_{j=i}^{n-1}\overline{k}_{j}^{n})^{2} &\leq &4\delta
E\sum_{j=i}^{n-1}[(u_{j})^{2}+g(t_{j},0,0)^{2}+\mu ^{2}(E[\bar{y}_{j+1}^{n}|%
\mathcal{G}_{j}^{n}])^{2}+\mu ^{2}(\bar{z}_{j}^{n})^{2}],  \nonumber
\end{eqnarray}
Set $\alpha =32\mu ^{2}$ in (\ref{est-imp2b3}), it follows
\begin{eqnarray*}
&&E|\bar{y}_{i}^{n}|^{2}+\frac{1}{4}\sum_{j=i}^{n-1}E|\bar{z}%
_{j}^{n}|^{2}\delta \\
&\leq &E|\xi ^{n}|^{2}+(\delta +\frac{1}{4\mu ^{2}}\delta +3\delta
^{2})E\sum_{j=i}^{n-1}[|g(t_{j},0,0)|^{2}]+32\mu
^{2}E[\sup_{j}((L_{j}^{n})^{+})^{2}+\sup_{j}((U_{j}^{n})^{+})^{2}] \\
&&+\delta (\frac{5}{4}+2\mu +4\mu ^{2}+3\mu ^{2}\delta )E\sum_{j=i}^{n-1}|%
\overline{y}_{j+1}^{n}|^{2}+\frac{1}{8\mu ^{2}}\delta
E\sum_{j=i}^{n-1}[(u_{i}^{l})^{2}+(u_{i}^{u})^{2}].
\end{eqnarray*}
Notice that $3\mu ^{2}\delta <1$, so $3\mu ^{2}\delta ^{2}<\delta $.
Then by applying the discrete Gronwall inequality in Lemma
\ref{dis-gro}, and the estimation of $\overline{a}_{j}^{n}$ and $\overline{k}%
_{j}^{n}$ follows from (\ref{est-2b}) , we get
\[
\sup_{j}E[\left| \overline{y}_{j}^{n}\right|
^{2}]+E[\sum_{j=0}^{n-1}\left|
\overline{z}_{j}^{n}\right| ^{2}\delta +\left| \sum_{j=0}^{n-1}\overline{k}%
_{j}^{n}\right| ^{2}+\left|
\sum_{j=0}^{n-1}\overline{a}_{j}^{n}\right| ^{2}]\leq c\;C_{\xi
^{n},g,L^{n},U^{n}}.
\]
We reconsider (\ref{sqr-exp}), as before take sum and $\sup_j$, then
take expectation, using Burkholder-Davis-Gundy inequality for
martingale part, with similar techniques, we get
\[
E[\sup_{j}\left| \overline{y}_{j}^{n}\right| ^{2}]\leq C_{\xi
^{n},g,L^{n},U^{n}}+C_{\mu}E[\sum_{j=0}^{n-1}\left|
\overline{y}_{j}^{n}\right| ^{2}\delta \leq E[\sup_{j}\left|
\overline{y}_{j}^{n}\right| ^{2}]\leq C_{\xi
^{n},g,L^{n},U^{n}}+C_{\mu}T\sup_j E[\left|
\overline{y}_{j}^{n}\right| ^{2}],
\]
which implies final result. $\square $

Then our convergence result for the explicit reflected scheme is

\begin{theorem}
Under the same assumptions as in Theorem \ref{conv-r2b}, when $n\rightarrow +\infty ,$%
\begin{equation}
E[\sup_{t}|\overline{Y}_{t}^{n}-Y_{t}|^{2}]+E\int_{0}^{T}|\overline{Z}%
_{t}^{n}-Z_{t}|^{2}dt\rightarrow 0.  \label{conv-exp-R2b}
\end{equation}
And $\overline{A}_{t}^{n}-\overline{K}_{t}^{n}\rightarrow A_{t}-K_{t}$ in $%
\mathbf{L}^{2}(\mathcal{F}_{t})$, for $0\leq t\leq T$.
\end{theorem}

\textbf{Proof.} Thanks to Theorem \ref{conv-r2b}, it is sufficient
to prove that
as $n\rightarrow +\infty ,$%
\begin{equation}
E[\sup_{j}|\overline{y}_{j}^{n}-y_{j}^{n}|^{2}]+E\sum_{j=0}^{n-1}|\overline{z%
}_{j}^{n}-z_{j}^{n}|^{2}\delta \rightarrow 0.  \label{conv-exp-R2}
\end{equation}
Since
\begin{eqnarray}
y_{j}^{n} &=&y_{j+1}^{n}+g(t_{j},y_{j}^{n},z_{j}^{n})%
\delta +a_{j}^{n}-k_{j}^{n}-z_{j}^n\sqrt{\delta}\varepsilon_{j+1}^n, \label{diff-ex-imp}\\
\bar{y}_{j}^{n} &=&E[\bar{y}_{j+1}^{n}|\mathcal{G}_{j}^{n}]+g(t_{j},E[\bar{y}%
_{j+1}^{n}|\mathcal{G}_{j}^{n}],\bar{z}_{j}^{n})\delta +\overline{a}_{j}^{n}-%
\overline{k}_{j}^{n}-\bar{z}_{j}^n\sqrt{\delta}\varepsilon_{j+1}^n,\nonumber
\end{eqnarray}
we get
\begin{eqnarray*}
&&E\left| y_{j}^{n}-\overline{y}_{j}^{n}\right| ^{2} \\
&=&E\left| y_{j+1}^{n}-\overline{y}_{j+1}^{n}\right| ^{2}-\delta E\left|
z_{j}^{n}-\overline{z}_{j}^{n}\right| ^{2}+2\delta E[(y_{j}^{n}-\overline{y}%
_{j}^{n})(g(t_{j},y_{j}^{n},z_{j}^{n})-g(t_{j},E[\overline{y}_{j+1}^{n}|%
\mathcal{G}_{j}^{n}],\overline{z}_{j}^{n}))] \\
&&-E[\delta (g(t_{j},y_{j}^{n},z_{j}^{n})-g(t_{j},E[\overline{y}_{j+1}^{n}|%
\mathcal{G}_{j}^{n}],\overline{z}_{j}^{n}))+(a_{j}^{n}-\overline{a}%
_{j}^{n})-(k_{j}^{n}-\overline{k}_{j}^{n})]^{2} \\
&&+2E[(y_{j}^{n}-\overline{y}_{j}^{n})(a_{j}^{n}-\overline{a}%
_{j}^{n})]-2E[(y_{j}^{n}-\overline{y}_{j}^{n})(k_{j}^{n}-\overline{k}%
_{j}^{n})] \\
&\leq &E\left| y_{j+1}^{n}-\overline{y}_{j+1}^{n}\right| ^{2}-\delta E\left|
z_{j}^{n}-\overline{z}_{j}^{n}\right| ^{2}+2\delta E[(y_{j}^{n}-\overline{y}%
_{j}^{n})(g(t_{j},y_{j}^{n},z_{j}^{n})-g(t_{j},E[\overline{y}_{j+1}^{n}|%
\mathcal{G}_{j}^{n}],\overline{z}_{j}^{n}))]
\end{eqnarray*}
in view of
\begin{eqnarray*}
(y_{j}^{n}-\overline{y}_{j}^{n})(a_{j}^{n}-\overline{a}_{j}^{n})
&=&(y_{j}^{n}-L_{j}^{n})a_{j}^{n}+(\overline{y}_{j}^{n}-L_{j}^{n})(\overline{%
a}_{j}^{n}) \\
&&-(\overline{y}_{j}^{n}-L_{j}^{n})a_{j}^{n}-(y_{j}^{n}-L_{j}^{n})(\overline{%
a}_{j}^{n}) \\
&\leq &0 \\
(y_{j}^{n}-\overline{y}_{j}^{n})(k_{j}^{n}-\overline{k}_{j}^{n})
&=&(y_{j}^{n}-U_{j}^{n})k_{j}^{n}+(\overline{y}_{j}^{n}-U_{j}^{n})\overline{k%
}_{j}^{n} \\
&&-(y_{j}^{n}-U_{j}^{n})\overline{k}_{j}^{n}-(\overline{y}%
_{j}^{n}-U_{j}^{n})k_{j}^{n} \\
&\geq &0.
\end{eqnarray*}
We take sum over $j$ from $i$ to $n-1$, with $\xi ^{n}-\overline{\xi }^{n}=0$%
, then we get
\begin{eqnarray*}
E\left| y_{j}^{n}-\overline{y}_{j}^{n}\right| ^{2}+\delta
\sum_{j=i}^{n-1}E\left| z_{j}^{n}-\overline{z}_{j}^{n}\right| ^{2} &\leq
&2\delta \sum_{j=i}^{n-1}E[(y_{j}^{n}-\overline{y}%
_{j}^{n})(g(t_{j},y_{j}^{n},z_{j}^{n})-g(t_{j},E[\overline{y}_{j+1}^{n}|%
\mathcal{G}_{j}^{n}],\overline{z}_{j}^{n}))] \\
&\leq &2\mu ^{2}\delta E\sum_{j=i}^{n-1}\left| y_{j}^{n}-\overline{y}%
_{j}^{n}\right| ^{2}+\frac{\delta }{2}\sum_{j=i}^{n-1}E\left| z_{j}^{n}-%
\overline{z}_{j}^{n}\right| ^{2} \\
&&+2\mu \delta E\sum_{j=i}^{n-1}\left| y_{j}^{n}-\overline{y}_{j}^{n}\right|
\cdot \left| y_{j}^{n}-E[\overline{y}_{j+1}^{n}|\mathcal{G}_{j}^{n}]\right| .
\end{eqnarray*}
Since $\bar{y}_{j}^{n}-E[\bar{y}_{j+1}^{n}|\mathcal{G}_{j}^{n}]=g(t_{j},E[%
\bar{y}_{j+1}^{n}|\mathcal{G}_{j}^{n}],\bar{z}_{j}^{n})\delta +\overline{a}%
_{j}^{n}-\overline{k}_{j}^{n}$, we have
\begin{eqnarray*}
&&2\mu \delta E[\left| y_{j}^{n}-\overline{y}_{j}^{n}\right| \cdot \left|
y_{j}^{n}-E[\overline{y}_{j+1}^{n}|\mathcal{G}_{j}^{n}]\right| ] \\
&=&2\mu \delta E[\left| y_{j}^{n}-\overline{y}_{j}^{n}\right| \cdot \left|
y_{j}^{n}-\overline{y}_{j}^{n}+g(t_{j},E[\overline{y}_{j+1}^{n}|\mathcal{G}%
_{j}^{n}]),z_{j}^{n})\delta +\overline{a}_{j}^{n}-\overline{k}%
_{j}^{n}\right| ] \\
&\leq &(2\mu +1)\delta E[\left| y_{j}^{n}-\overline{y}_{j}^{n}\right|
^{2}]+\mu ^{2}\delta E[3\delta ^{2}(\left| g(t_{j},0,0)\right| ^{2}+\mu
^{2}\left| \overline{y}_{j+1}^{n}\right| ^{2}+\mu ^{2}\left|
z_{j}^{n}\right| ^{2}+(\overline{a}_{j}^{n})^{2}+(\overline{k}_{j}^{n})^{2}].
\end{eqnarray*}
Then by Lemma \ref{est-expRBSDE2d}, we obtain
\begin{equation}
E\left| y_{j}^{n}-\overline{y}_{j}^{n}\right| ^{2}+\frac{\delta }{2}%
\sum_{j=i}^{n-1}E\left| z_{j}^{n}-\overline{z}_{j}^{n}\right| ^{2}\leq (2\mu
^{2}+2\mu +1)\delta \sum_{j=i}^{n-1}E\left| y_{j}^{n}-\overline{y}%
_{j}^{n}\right| ^{2}+\delta C_{\xi ^{n},g,L^{n},U^{n}}.  \label{est-exp-d2}
\end{equation}
By the discrete Gronwall inequality in Lemma \ref{dis-gro}, we get
\[
\sup_{j\leq n}E\left| y_{j}^{n}-\overline{y}_{j}^{n}\right| ^{2}\leq C\delta
^{2}e^{(2\mu +2\mu ^{2}+1)T}.
\]
 With (\ref{est-exp-d2}), it follows $E[\delta%
\sum_{j=0}^{n-1}E\left| z_{j}^{n}-\overline{z}_{j}^{n}\right|
^{2}]\leq C\delta^2$.
 Then we reconsider (\ref{diff-ex-imp}), this time we take
expectation after taking square, sum and $\sup$ over $j$. Using
Burkholder-Davis-Gundy inequality for martingale parts and similar
tachniques, it follows that
\[
E \sup_{j\leq n}\left| y_{j}^{n}-\overline{y}_{j}^{n}\right|
^{2}\leq CE[\sum_{j=0}^{n-1}E\left|
y_{j}^{n}-\overline{y}_{j}^{n}\right| ^{2}\delta]\leq CT \sup_{j\leq
n}E\left| y_{j}^{n}-\overline{y}_{j}^{n}\right| ^{2},
\]
 which implies (\ref{conv-exp-R2}).

For the convergence of $(\overline{A}^{n},\overline{K}^{n})$, we consider
\begin{eqnarray*}
\overline{A}_{t}^{n}-\overline{K}_{t}^{n} &=&\overline{Y}_{0}^{n}-\overline{Y%
}_{t}^{n}-\int_{0}^{t}g(s,\overline{Y}_{s}^{n},\overline{Z}%
_{s}^{n})ds+\int_{0}^{t}\overline{Z}_{s}^{n}dB_{s}^{n}, \\
A_{t}^{n}-K_{t}^{n}
&=&Y_{0}^{n}-Y_{t}^{n}-\int_{0}^{t}g(s,Y_{s}^{n},Z_{s}^{n})ds+%
\int_{0}^{t}Z_{s}^{n}dB_{s}^{n},
\end{eqnarray*}
then the convergence results follow easily from the convergence of
$A^n$, (\ref{conv-exp-R2b}) and the Lipschitz condition of $g$.
$\square $

\section{Simulations of Reflected BSDE with two barriers}

For computational convenience, we consider the case when $T=1$. The
calculation begins from $y_{n}^{n}=\xi ^{n}$ and proceeds backward
to solve $(y_{j}^{n},z_{j}^{n},a_{j}^{n},k_{j}^{n})$, for
$j=n-1,\cdots 1,0$. Due to the amount of computation, we consider a
very simple case: $\xi
=\Phi (B_{1})$, $L_{t}=\psi _{1}(t,B(t))$, $U_{t}=\psi _{2}(t,B(t))$, where $%
\Phi $, $\psi _{1}$ and $\psi _{2}$ are real analytic functions defined on $%
\mathbb{R }$ and $[0,1]\times \mathbb{R}$ respectively. As mentioned
in the introduction, we have developed a Matlab toolbox for
calculating and simulating solutions of reflected BSDEs with two
barriers which has a well-designed interface. This toolbox can be
downloaded from
http://159.226.47.50:8080/iam/xumingyu/English/index.jsp, with
clicking 'Preprint' on the left side.

We take the following example: $g(y,z)=-5\left| y+z\right| -1$,
$\Phi (x)=\left|
x\right| $, $\Psi _{1}(t,x)=-3(x-2)^{2}+3$, $\Psi _{2}(t,x)=(x+1)^{2}+3(t-1)$%
, and $n=400$.

In Figure 1, we can see both the global situation of the solution surface of $%
y^{n}$ and its partial situation i.e. trajectory. In the upper portion of
Figure 1, it is in 3-dimensional. The lower surface shows the barrier $L$%
, as well the upper one is for the barrier $U$. The solution $y^{n}$
is in the middle of them. Then we generate one trajectory of the
discrete Brownian motion $(B_{j}^{n})_{0\leq j\leq n}$, which is
drawn on the horizontal plane. The value of $y_{j}^{n}$ with respect
to this Brownian sample is showed on the solution surface. The
remainder of the figure shows respectively the
trajectory of the force $A_{j}^{n}=\sum_{i=0}^{j}a_{i}^{n}$ and $%
K_{j}^{n}=\sum_{i=0}^{j}k_{i}^{n}.$

\begin{center}
\centering\includegraphics[totalheight=110mm]{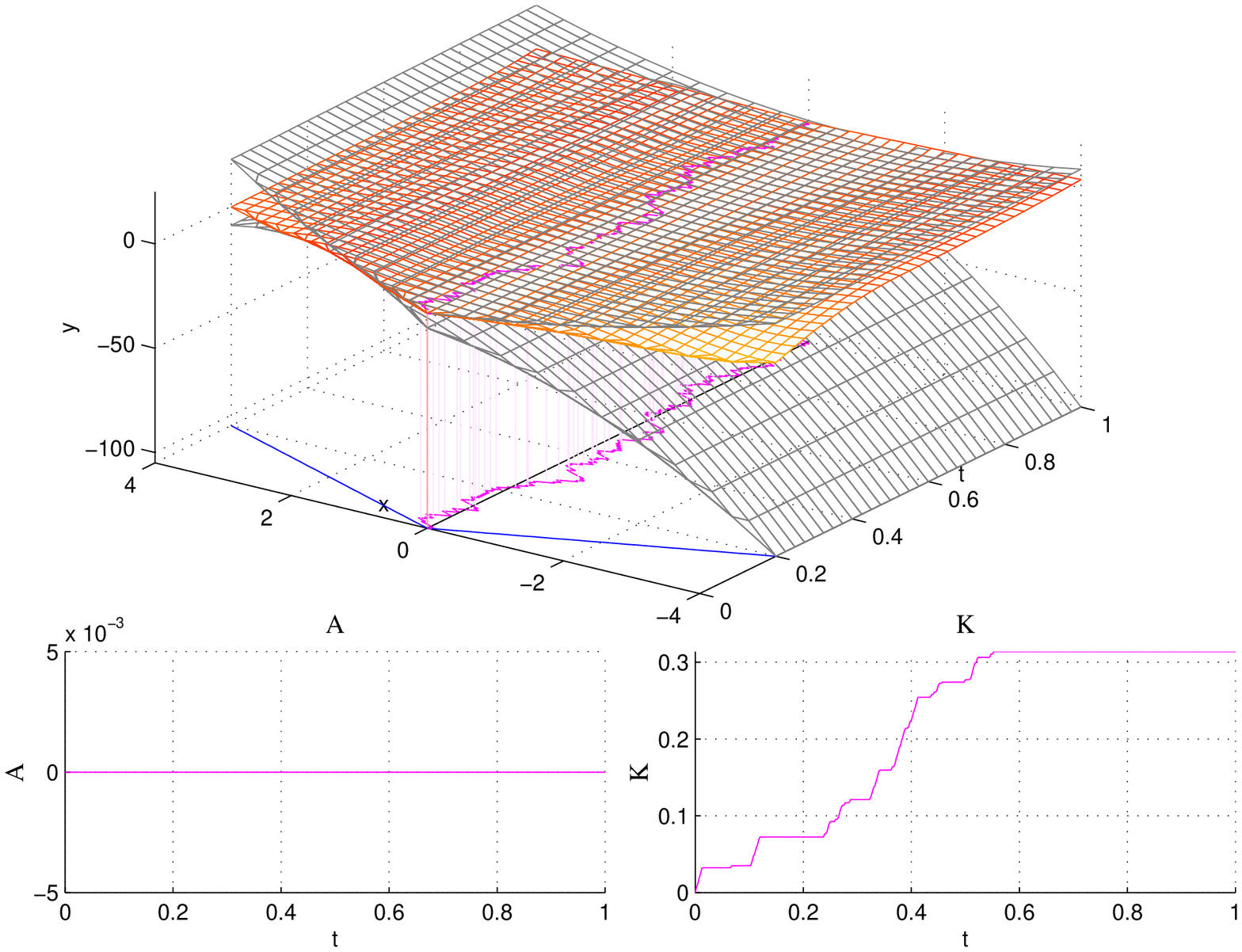} \\[0pt]
Figure 1: A solution surface of reflected BSDE with two barriers
\end{center}
The lower graphs shows clearly that $%
A^{n} $ (respective $K^{n}$) acts only if $y^{n}$ touches the lower barrier $%
L^{n}$, i.e. on the set $\{y^{n}=L^{n}\}$ (respective the upper
barrier $U$, i.e. on the set $\{y^{n}=U^{n}\}$ ), and they never act
at the same time.

 In the upper portion we can see that there is an area,
named Area I, (resp. Area II) where the solution surface and the
lower barrier surface (resp. the solution surface and the upper
barrier surface) stick together. When the trajectory of solution
$y_{j}^{n}$ goes into Area I (resp. Area II), the force $A_{j}^{n}$
(resp. $K_{j}^{n}$) will push $y_{j}^{n}$ upward (resp. downward).
Indeed, if we don't have the barriers here, $y_{j}^{n}$ intends
going up or down to cross the reflecting barrier $L_{j}^{n}$ and $U_{j}^{n}$%
, so to keep $y_{j}^{n}$ staying between $L_{j}^{n}$ and
$U_{j}^{n}$, the action of forces $A_{j}^{n}$ and $K_{j}^{n}$ are
necessary. In Figure 1, the increasing process $A_{j}^{n}$ keeps
zero, while $K_{j}^{n}$ increases from the beginning.
Correspondingly in the beginning $y_{j}^{n}$ goes into Area II, but
always stay out of Area I. Since Area I and Area II are totally
disjoint, so $A_{j}^{n}$ and $K_{j}^{n}$ never increase at same
time.

About this point, let us have a look at Figure 2. This figure shows
a group of
3-dimensional dynamic trajectories $(t_{j},B_{j}^{n},Y_{j}^{n})$ and $%
(t_{j},B_{j}^{n},Z_{j}^{n})$, simultaneously, of 2-dimensional
trajectories of $(t_{j},Y_{j}^{n})$ and $(t_{j},Z_{j}^{n})$. For the
other sub-figures, the upper-right one is for the trajectories
$A_{j}^{n}$, and while the lower-left one is for $K_{j}^{n}$, then
comparing these two sub-figures, as in Figure 1, $\{a_{j}^{n}\neq
0\}$ and $\{k_{j}^{n}\neq 0\}$ are disjoint, but the converse is not
true.

\begin{center}
\centering\includegraphics[totalheight=120mm]{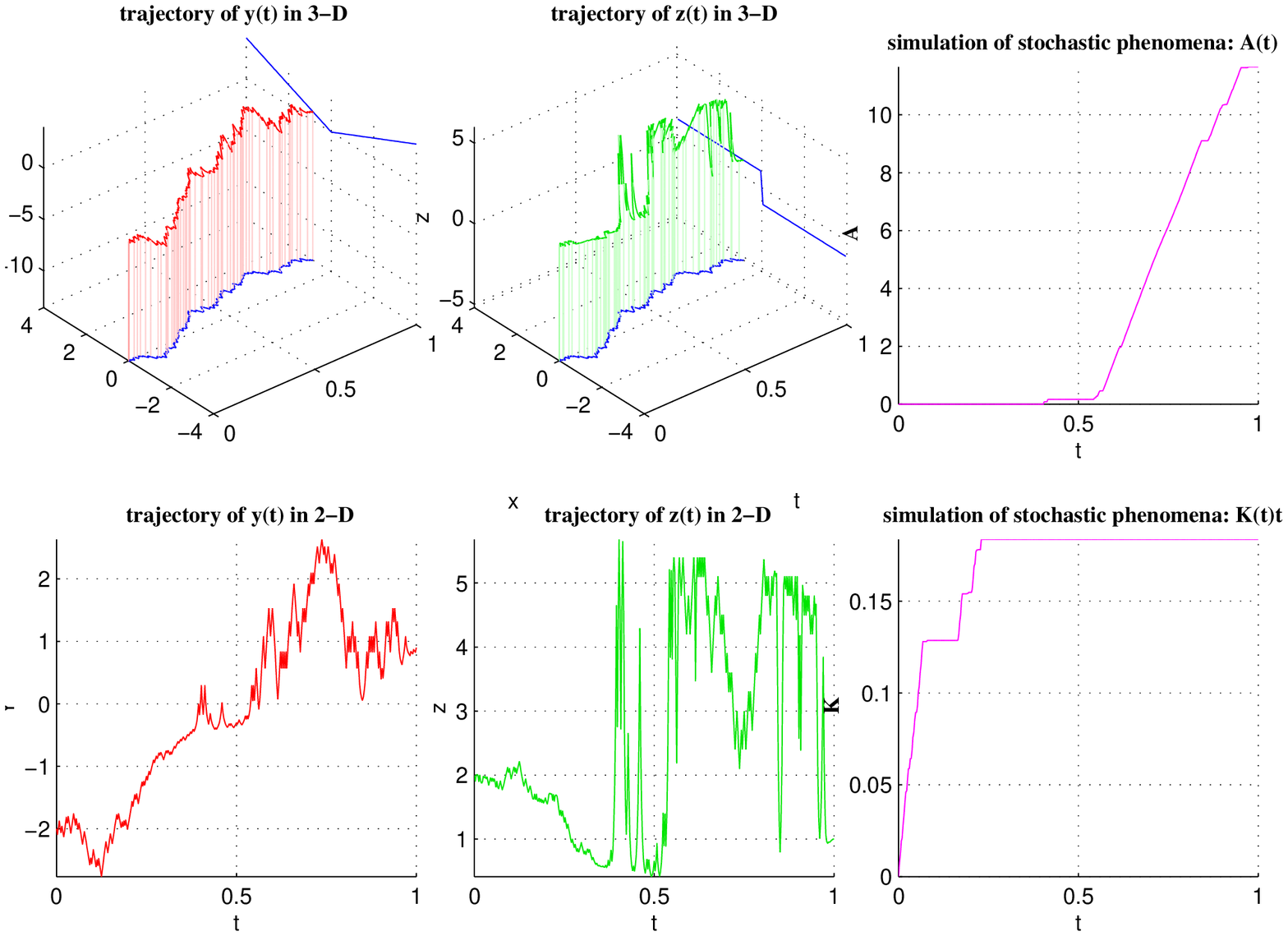}\\[0pt]
Figure 2: The trajectories of solutions of (\ref{RBSDE2b1})
\end{center}

Now we present some numerical results using the explicit reflected
scheme and the implicit-explicit penalization scheme, respectively,
with different discretization. Consider the
parameters: $g(y,z)=-5\left| y+z\right| -1$, $\Phi (x)=\left| x\right| $, $%
\Psi _{1}(t,x)=-3(x-2)^{2}+3$, $\Psi _{2}(t,x)=(x+1)^{2}+3t-2.5$:
\begin{eqnarray*}
&&
\begin{array}{lll}
n=400, & \mbox{reflected explicit scheme:} & y_{0}^{n}=-1.7312 \\[0.3cm]
& \mbox{penalization scheme:} &
\begin{tabular}{|c|c|c|c|c|}
\hline
$p$ & $20$ & $200$ & $2000$ & $2\times 10^{4}$ \\ \hline
$y_{0}^{p,n}$ & $-1.8346$ & $-1.7476$ & $-1.7329$ & $-1.7314$ \\ \hline
\end{tabular}
\end{array}
\\
&&
\begin{array}{lll}
n=1000, & \mbox{reflected explicit scheme:} & y_{0}^{n}=-1.7142 \\[0.3cm]
& \mbox{penalization scheme:} &
\begin{tabular}{|c|c|c|c|c|}
\hline
$p$ & $20$ & $200$ & $2000$ & $2\times 10^{4}$ \\ \hline
$y_{0}^{p,n}$ & $-1.8177$ & $-1.7306$ & $-1.7161$ & $-1.7144$ \\ \hline
\end{tabular}
\end{array}
\\
&&
\begin{array}{lll}
n=2000, & \mbox{reflected explicit scheme:} & y_{0}^{n}=-1.7084 \\[0.3cm]
& \mbox{penalization scheme:} &
\begin{tabular}{|c|c|c|c|c|}
\hline
$p$ & $20$ & $200$ & $2000$ & $2\times 10^{4}$ \\ \hline
$y_{0}^{p,n}$ & $-1.8124$ & $-1.7250$ & $-1.7103$ & $-1.7068$ \\ \hline
\end{tabular}
\end{array}
\\
&&
\begin{array}{lll}
n=4000, & \mbox{reflected explicit scheme:} & y_{0}^{n}=-1.7055 \\[0.3cm]
& \mbox{penalization scheme:} &
\begin{tabular}{|c|c|c|c|c|}
\hline
$p$ & $20$ & $200$ & $2000$ & $2\times 10^{4}$ \\ \hline
$y_{0}^{p,n}$ & $-1.8096$ & $-1.7222$ & $-1.7074$ & $-1.7057$ \\ \hline
\end{tabular}
\end{array}
\end{eqnarray*}

From this form, first we can see that as the penalization parameter
$p$ increases, the penalization solution $y_{0}^{p,n}$ tends
increasingly to the reflected solution $y_{0}^{n}$. Second, as the
discretaization parameter $n$ increases, the differences of
$y_{0}^{n}$ with different $n$ become smaller as well as that of
$y_{0}^{p,n}$.

\section{Appendix: The proof of Theorem \ref{p2b}}

To complete the paper, here we give the proof of Theorem \ref{p2b}.

\noindent\textbf{Proof of Theorem \ref{p2b}.} (a) is the main result
in \cite{LSM}. So we omit its proof.

Now we consider (b). The convergence of $(Y_{t}^{p},Z_{t}^{p})$ is a
direct consequence of \cite{PX2003}. For the convergence speed, the
proof is a combination of results in \cite{LSM} and \cite{PX2003}.
From \cite{LSM}, we
 know that for (\ref{PBSDE2b-g}), as $m\rightarrow \infty $, $%
\widehat{Y}_{t}^{m,p}\nearrow \underline{Y}_{t}^{p}$ in $\mathbf{S}^{2}(0,T)$, $%
\widehat{Z}_{t}^{m,p}\rightarrow \underline{Z}_{t}^{p}$ in $\mathbf{L}_{\mathcal{F}%
}^{2}(0,T)$, $\widehat{A}_{t}^{m,p}\rightarrow \underline{A}_{t}^{p}$ in $\mathbf{S}%
^{2}(0,T)$, where $(\underline{Y}_{t}^{p},\underline{Z}_{t}^{p},\underline{A}%
_{t}^{p})$ is a solution of the following reflected BSDE with one lower barrier $L$%
\begin{eqnarray}
-d\underline{Y}_{t}^{p} &=&g(t,\underline{Y}_{t}^{p},\underline{Z}%
_{t}^{p})dt+d\underline{A}_{t}^{p}-p(\underline{Y}_{t}^{p}-U_{t})^{+}dt-%
\underline{Z}_{t}^{p}dB_{t},\underline{Y}_{T}^{p}=\xi ,  \label{PBSDE2b-l} \\
\underline{Y}_{t}^{p} &\geq &L_{t},\int_{0}^{T}(\underline{Y}_{t}^{p}-L_{t})d%
\underline{A}_{t}^{p}=0.  \nonumber
\end{eqnarray}
Set $\underline{K}_{t}^{p}=\int_{0}^{t}p(\underline{Y}_{s}^{p}-U_{s})^{+}ds$%
. Then letting $p\rightarrow \infty $, it follows that $\underline{Y}%
_{t}^{p}\searrow Y_{t}$ in $\mathbf{S}^{2}(0,T)$, $\underline{Z}%
_{t}^{p}\rightarrow Z_{t}$ in $\mathbf{L}_{\mathcal{F}}^{2}(0,T)$. By
comparison theorem, $d\underline{A}_{t}^{p}$ is increasing in $p$. So $%
\underline{A}_{T}^{p}\nearrow A_{T}$, and $0\leq \sup_{t}[\underline{A}%
_{t}^{p+1}-\underline{A}_{t}^{p}]\leq \underline{A}_{T}^{p+1}-\underline{A}%
_{T}^{p}$. It follows that $\underline{A}_{t}^{p}\rightarrow A_{t}$ in $%
\mathbf{S}^{2}(0,T)$. Then with Lipschitz condition of $g$ and convergence
results, we get $\underline{K}_{t}^{p}\rightarrow K_{t}$ in $\mathbf{S}%
^{2}(0,T)$. Moreover from Lemma 4 in \cite{LSM}, we know that there
exists a constant $C$ depending on $\xi $, $g(t,0,0)$, $\mu $, $L$
and $U$, such that
\[
E[\sup_{0\leq t\leq T}|\underline{Y}_{t}^{p}-Y_{t}|^{2}+\int_{0}^{T}(|%
\underline{Z}_{t}^{p}-Z_{t}|^{2})dt\leq \frac{C}{\sqrt{p}}.
\]
Similarly for (\ref{PBSDE2b-g}), first letting $p\rightarrow \infty $, we
get $\widehat{Y}_{t}^{m,p}\searrow \overline{Y}_{t}^{m}$ in $\mathbf{S}^{2}(0,T)$, $%
\widehat{Z}_{t}^{m,p}\rightarrow \overline{Z}_{t}^{m}$ in $\mathbf{L}_{\mathcal{F}%
}^{2}(0,T)$, $\widehat{K}_{t}^{m,p}\rightarrow \overline{K}_{t}^{m}$ in $\mathbf{S}%
^{2}(0,T)$, where $(\overline{Y}_{t}^{m},\overline{Z}_{t}^{m},\overline{K}%
_{t}^{m})$ is a solution of the following reflected BSDE with one upper barrier $U$%
\begin{eqnarray}
-d\overline{Y}_{t}^{m} &=&g(t,\overline{Y}_{t}^{m},\overline{Z}%
_{t}^{m})dt+m(L_{t}-\overline{Y}_{t}^{m})^{+}dt-d\overline{K}_{t}^{m}-%
\overline{Z}_{t}^{m}dB_{t},\overline{Y}_{T}^{m}=\xi ,  \label{PBSDE2b-u} \\
\overline{Y}_{t}^{m} &\leq &U_{t},\int_{0}^{T}(\overline{Y}_{t}^{p}-U_{t})d%
\overline{K}_{t}^{m}=0.  \nonumber
\end{eqnarray}
In the same way, as $m\rightarrow \infty $, $\overline{Y}_{t}^{m}\nearrow
Y_{t}$ in $\mathbf{S}^{2}(0,T)$, $\overline{Z}_{t}^{m}\rightarrow Z_{t}$ in $%
\mathbf{L}_{\mathcal{F}}^{2}(0,T)$, and $(\overline{A}_{t}^{m},\overline{K}%
_{t}^{m})\rightarrow (A_{t},K_{t})$ in $(\mathbf{S}^{2}(0,T))^{2}$, where $%
\overline{A}_{t}^{m}=\int_{0}^{t}m(L_{s}-\overline{Y}_{s}^{m})^{+}ds$. Also
there exists a constant $C$ depending on $\xi $, $g(t,0,0)$, $\mu $, $L$ and
$U$, such that
\[
E[\sup_{0\leq t\leq T}|\overline{Y}_{t}^{m}-Y_{t}|^{2}+\int_{0}^{T}(|%
\overline{Z}_{t}^{m}-Z_{t}|^{2})dt\leq \frac{C}{\sqrt{m}}.
\]

Applying comparison theorem to (\ref{PBSDE2b}) and
(\ref{PBSDE2b-l}), (\ref
{PBSDE2b}) and (\ref{PBSDE2b-u})(let $m=p$), we have $\overline{Y}%
_{t}^{p}\leq Y_{t}^{p}\leq \underline{Y}_{t}^{p}$. Then we have
\[
E[\sup_{0\leq t\leq T}|Y_{t}^{p}-Y_{t}|^{2}]\leq \frac{C}{\sqrt{p}},
\] for some constant $C$. To get the estimate results for $Z^{p}$, we apply
It\^{o} formula to $\left| Y_{t}^{p}-Y_{t}\right| ^{2}$, and get
\begin{eqnarray*}
&&E|Y_{0}^{p}-Y_{0}|^{2}+\frac{1}{2}E\int_{0}^{T}|Z_{s}^{p}-Z_{s}|^{2}ds \\
&=&(\mu +2\mu ^{2})E\int_{0}^{T}\left| Y_{s}^{p}-Y_{s}\right|
^{2}ds+2E\int_{0}^{T}(Y_{s}^{p}-Y_{s})dA_{s}^{p}-2E%
\int_{0}^{T}(Y_{s}^{p}-Y_{s})dA_{s} \\
&&-2E\int_{0}^{T}(Y_{s}^{p}-Y_{s})dK_{s}^{p}+2E%
\int_{0}^{T}(Y_{s}^{p}-Y_{s})dK_{s}.
\end{eqnarray*}
Since
\begin{eqnarray*}
2E\int_{0}^{T}(Y_{s}^{p}-Y_{s})dA_{s}^{p}
&=&2E\int_{0}^{T}(Y_{s}^{p}-L_{s})dA_{s}^{p}-2E%
\int_{0}^{T}(Y_{s}-L_{s})dA_{s}^{p} \\
&\leq &2pE\int_{0}^{T}(Y_{s}^{p}-L_{s})(Y_{s}^{p}-L_{s})^{-}ds\leq 0
\end{eqnarray*}
and $2E\int_{0}^{T}(Y_{s}^{p}-Y_{s})dK_{s}^{p}\geq
2pE\int_{0}^{T}(Y_{s}^{p}-U_{s})(Y_{s}^{p}-U_{s})^{+}ds\geq 0$, we have
\[
E\int_{0}^{T}|Z_{s}^{p}-Z_{s}|^{2}ds\leq \frac{C}{\sqrt{p}},
\]
in view of the estimation of $A$ and $K$ and the convergence of
$Y^{p}$.

Now we consider the convergence of $A^{p}$ and $K^{p}$. Since
\begin{eqnarray*}
A_{t}-K_{t}
&=&Y_{0}-Y_{t}-\int_{0}^{t}g(s,Y_{s},Z_{s})ds+\int_{0}^{t}Z_{s}dB_{s}, \\
A_{t}^{p}-K_{t}^{p}
&=&Y_{0}^{p}-Y_{t}^{p}-\int_{0}^{t}g(s,Y_{s}^{p},Z_{s}^{p})ds+%
\int_{0}^{t}Z_{s}^{p}dB_{s},
\end{eqnarray*}
from the Lipschitz condition of $g$ and the convergence results of $Y^{p}$ and $%
Z^{p} $, we have immediately
\begin{eqnarray*}
&&E[\sup_{0\leq t\leq T}[(A_{t}-K_{t})-(A_{t}^{p}-K_{t}^{p})]^{2}] \\
&\leq &8E[\sup_{0\leq t\leq T}|Y_{t}^{p}-Y_{t}|^{2}+4\mu
\int_{0}^{T}|Y_{s}^{p}-Y_{s}|^{2}ds+C\int_{0}^{T}|Z_{s}^{p}-Z_{s}|^{2}ds]%
\leq \frac{C}{\sqrt{p}}.
\end{eqnarray*}
Meanwhile we know $E[(A_{T}^{p})^{2}+(K_{T}^{p})^{2}]<\infty $, so $A^{p}$
and $K^{p}$ admits weak limit $\widetilde{A}$ and $\widetilde{K}$ in $%
\mathbf{S}^{2}(0,T)$ respectively. By the comparison results of $\overline{Y}%
_{t}^{p}$, $Y_{t}^{p}$ and $\underline{Y}_{t}^{p}$, we get
\begin{eqnarray*}
dA_{t}^{p} &=&p(Y_{t}^{p}-L_{t})^{-}dt\leq p(\overline{Y}%
_{t}^{p}-L_{t})^{-}dt=d\overline{A}_{t}^{p}, \\
dK_{t}^{p} &=&p(Y_{t}^{p}-U_{t})^{+}dt\geq p(\underline{Y}%
_{t}^{p}-U_{t})^{+}dt=d\underline{K}_{t}^{p}.
\end{eqnarray*}
So $d\widetilde{A}_{t}\leq dA_{t}$ and $d\widetilde{K}_{t}\geq dK_{t}$, it
follows that $d\widetilde{A}_{t}-d\widetilde{K}_{t}\leq dA_{t}-dK_{t}$. On
the other hand, the limit of $Y^{p}$ is $Y$, so $d\widetilde{A}_{t}-d%
\widetilde{K}_{t}=dA_{t}-dK_{t}$. Then there must be $d\widetilde{A}%
_{t}=dA_{t}$ and $d\widetilde{K}_{t}=dK_{t}$, which implies $\widetilde{A}%
_{t}=A_{t}$ and $\widetilde{K}_{t}=K_{t}$. $\square $

\end{document}